\documentclass[11pt]{amsart}
\usepackage{amssymb,amsmath,amsthm}
\newtheorem{theorem}{Theorem}[section]

\newtheorem{cor}[theorem]{Corollary}
\newtheorem{conj}[theorem]{Conjecture}
\newtheorem{prop}[theorem]{Proposition}

\theoremstyle{definition}

\theoremstyle{remark}

\numberwithin{equation}{section}

\newcommand{\abs}[1]{\lvert#1\rvert}
\newcommand{\spt}{\mbox{\rm spt}}
\newcommand{\SPT}{\mbox{\rm SPT}}
\newcommand{\Span}{\mbox{Span}}

\newcommand\Etwid{\overset {\text{\lower 3pt\hbox{$\sim$}}}E}
\newcommand\Ftwid{\overset {\text{\lower 3pt\hbox{$\sim$}}}F}
\newcommand\Qtwid{\overset {\text{\lower 3pt\hbox{$\sim$}}}Q}

\newcommand{\Uop}[2]{U_{0,#1}^{*}\left(#2\right)}

\DeclareSymbolFont{AMSb}{U}{msb}{m}{n}
\DeclareMathSymbol{\Z}{\mathalpha}{AMSb}{"5A}

\DeclareMathSymbol{\nmid}{\mathrel}{AMSb}{"2D}
\DeclareSymbolFont{AMSb}{U}{msb}{m}{n}
\DeclareMathSymbol{\C}{\mathalpha}{AMSb}{"43}
\DeclareMathSymbol{\F}{\mathalpha}{AMSb}{"46}
\DeclareMathSymbol{\N}{\mathalpha}{AMSb}{"4E}
\DeclareMathSymbol{\Q}{\mathalpha}{AMSb}{"51}
\DeclareMathSymbol{\R}{\mathalpha}{AMSb}{"52}
\DeclareMathSymbol{\Z}{\mathalpha}{AMSb}{"5A}

\newcommand{\Hup}{\mathcal{H}} 

\allowdisplaybreaks

\begin{document}
\newcommand{\beqs}{\begin{equation*}}
\newcommand{\eeqs}{\end{equation*}}
\newcommand{\beq}{\begin{equation}}
\newcommand{\eeq}{\end{equation}}
\newcommand\mylabel[1]{\label{#1}}
\newcommand\eqn[1]{(\ref{eq:#1})}
\newcommand\thm[1]{\ref{thm:#1}}
\newcommand\lem[1]{\ref{lem:#1}}
\newcommand\propo[1]{\ref{propo:#1}}
\newcommand\corol[1]{\ref{cor:#1}}
\newcommand\sect[1]{\ref{sec:#1}}
\newcommand\gausspoly[2]{\begin{bmatrix} #1 \\ #2\end{bmatrix}_q}
\newcommand\mytwid[1]{\overset {\text{\lower 3pt\hbox{$\sim$}}}#1}
\newcommand\leg[2]{\genfrac{(}{)}{}{}{#1}{#2}} 

\title[Congruences for the Smallest Parts Function and the Rank]{Congruences
for Andrews' Smallest Parts Partition \\
Function and New Congruences for Dyson's Rank}

\author{F.G. Garvan}
\address{Department of Mathematics, University of Florida, Gainesville,
Florida 32611-8105}
\email{frank@math.ufl.edu}          
\thanks{The author was supported in part by NSA Grant H98230-07-1-0011.
The original version of this paper was written July 1, 2007.
The main results of this paper were first presented at 
the Illinois Number Theory Fest,
May 17, 2007.
}



\date{June 11, 2008}%



\begin{abstract}
Let $\spt(n)$ denote the total number of appearances of smallest parts in the
partitions of $n$. Recently, Andrews showed how $\spt(n)$ is related to the 
second rank moment, and proved some surprising Ramanujan-type congruences 
mod $5$, $7$ and $13$.  We prove a generalization of these congruences
using known relations between rank and crank moments. 
We obtain explicit Ramanujan-type congruences for $\spt(n)$ mod $\ell$
for $\ell=11$, $17$, $19$, $29$, $31$ and $37$. 
Recently, Bringmann and Ono proved that Dyson's rank function
has infinitely many Ramanujan-type congruences. Their proof is
non-constructive and utilizes the theory of weak Maass forms.
We construct two explicit nontrivial examples mod $11$ using elementary 
congruences between rank moments and half-integer weight Hecke eigenforms.
\end{abstract}

\maketitle

\section{Introduction} \label{sec:intro}

Let $\spt(n)$ denote the number of smallest parts in the partitions of $n$.
Below is a list of the partitions of $4$ with their corresponding number 
of smallest parts. We see that $\spt(4)=10$.
$$
\begin{array}{ll}
4 & 1 \\
3+1 & 1\\
2+2 & 2\\
2+1+1 & 2\\
1+1+1+1 & 4
\end{array}
$$
In a recent paper Andrews \cite{Aspt} showed that $\spt(n)$ is
related to the second rank moment. The rank of a partition \cite{Dyson44}
is the
largest part minus the number of parts. The crank of a partition \cite{AG} is
the largest part if the partition has no ones, otherwise it is the
difference between the number of parts larger than the number of ones
and the number of ones. We let $N(m,n)$ denote the number of 
partitions of $n$ with rank $m$. For $n\ne1$ we let $M(m,n)$
denote the number of partitions of $n$ with crank $m$. For $n=1$
we define
$$
M(-1,1)=1, M(0,1)=-1, M(1,1)=1,\quad\mbox{and otherwise $M(m,1)=0$.}
$$
For $k$ even we define
\begin{align*}
N_k(n) &= \sum_{m} m^k N(m,n) \qquad\mbox{($k$th rank moment)},\\
M_k(n) &= \sum_{m} m^k M(m,n) \qquad\mbox{($k$th crank moment)}.
\end{align*}
Andrews proved that
\beq
\spt(n) = n\,p(n) - \tfrac12 N_2(n),
\mylabel{eq:AndSptId}
\eeq
where $p(n)$ is the number of partitions of $n$.
Dyson \cite{Dy} gave a combinatorial proof that
\beq
\tfrac12 M_2(n) = n\,p(n).
\mylabel{eq:DyM2}
\eeq
Hence we have
\beq
\spt(n) = \tfrac12(M_2(n) - N_2(n)).
\mylabel{eq:sptid}
\eeq
We make the following
\begin{conj}
$$
M_k(n) > N_k(n),
$$
for $k$ even, $k\ge2$ and $n\ge1$.
\end{conj}
The case $k=2$ follows from \eqn{sptid}.
We have checked the conjecture for $k\le 10$ and $n\le 500$.

In \cite{Aspt}, Andrews found some surprising congruences for 
$\spt(n)$
\begin{align}
\spt(5n + 4) &\equiv 0 \pmod{5},\mylabel{eq:spt5}\\
\spt(7n + 5) &\equiv 0 \pmod{7},\mylabel{eq:spt7}\\
\spt(13n + 6) &\equiv 0 \pmod{13}\mylabel{eq:spt13}.
\end{align}
These congruences are reminiscent of Ramanujan's partition
congruences
\begin{align}
p(5n + 4) &\equiv 0 \pmod{5},\mylabel{eq:p5}\\
p(7n + 5) &\equiv 0 \pmod{7},\mylabel{eq:p7}\\
p(11n + 6) &\equiv 0 \pmod{11}.\mylabel{eq:p11}
\end{align}
Andrews' proof of \eqn{spt5}, \eqn{spt7} depends solely on
\eqn{p5}, \eqn{p7} and known relations for the rank mod $5$ and $7$.
The proof of \eqn{spt13} is more difficult and depends on
relations mod $13$ for the rank due to O'Brien \cite{OB}.

In Section \sect{spt5713} we prove the following generalization 
of \eqn{spt5}--\eqn{spt13}:
For $t=5$, $7$ or $13$,     
\beqs
\spt(n) \equiv -2(n + \tfrac{t^2-1}{24})\,p(n) \pmod{t},
\eeqs
provided $1-24n$ is not a quadratic residue mod $t$.
The proof uses known 
relations between rank and crank moments \cite{AtG}. 

In Section \sect{modforms} we prove
a number of results on Hecke operators and congruences for modular forms
that will be needed in later sections.
In Section \sect{spt11} we study $\spt(n)\pmod{11}$.
By using a known relation between rank moments, crank moments
and the 23rd power of the eta-function we find that
\beqs
\sum_{n=0}^\infty \spt(11n+6) q^n\equiv 4 \prod_{n=1}^\infty (1-q^n)^{13}
\pmod{11}.
\eeqs
By using the fact that $\eta(24\tau)^{13}$ is a Hecke eigenform
we obtain the congruence
\beqs
\spt( 11\cdot19^4\cdot n + 22006) \equiv 0 \pmod{11}. 
\eeqs
In Section \sect{explicitrankcongs} we show that the 
generating functions for the rank moments $N_{2k}(11n+6) \pmod{11}$
can basically be written in terms of half-integral weight Hecke
eigenforms.
As a result we are able to find the following explicit congruences
for the rank
\begin{align}
N(r,11,5^4\cdot11\cdot19^4\cdot n + 4322599) &\equiv 0 \pmod{11},
\mylabel{eq:Intronicecong11a}\\
N(r,11,11^2\cdot19^4\cdot n + 172904) &\equiv 0 \pmod{11},
\mylabel{eq:Intronicecong11b}
\end{align}
for all $0\le r \le 10$.
Bringmann and Ono \cite{BringmannOno} conjectured and  
Bringmann \cite{Bringmann} proved that for each prime $\ell>3$ 
and $m$, $n\in \N$ there are infinitely many non-nested arithmetic
progressions $An+B$ such that
\beq
N(r,\ell^m,An+B) \equiv 0 \pmod{\ell^u},
\mylabel{eq:NrAnBcongs}
\eeq
for all $0 \le r \le \ell^{m}-1$. The congruences \eqn{Intronicecong11a} and
\eqn{Intronicecong11b} represent the first nontrivial explicit
examples of this result. The analogue of \eqn{NrAnBcongs}
for $N(r,t,An+B)$ when $(t,2\ell)=1$ had been proved earlier
by Bringmann and Ono \cite{BringmannOno}. The crank analogue was proved
by Mahlburg \cite{Mahlburg}. All these results generalize the analog
for $p(n)$ proved by Ono \cite{Ono}. It is also clear by \eqn{AndSptId} that
$\spt(n)$ can be written as an integer-linear combination of the rank
functions $N(r,\ell^u,n)$ mod $\ell^u$. Thus Bringmann's result \eqn{NrAnBcongs}
also implies that there are infinitely many non-nested arithmetic
progressions $An+B$ such that
\beq
\spt(An+B) \equiv 0 \pmod{\ell^u}.
\mylabel{eq:SptAnBcongs}
\eeq
In Section \sect{explicitsptcongs}
we describe algorithms for computing congruences for the generating
functions of $\spt(\ell n + \beta_\ell) \pmod{\ell}$ when 
$24\beta_\ell\equiv1\pmod{\ell}$, and Ramanujan-type congruences
for $\spt(n)$. For certain small primes we find that an appropriate form for
the generating function of $\spt(\ell n + \beta_\ell)$ is congruent mod $\ell$
to a half-integer weight Hecke eigenform. This leads to explicit 
examples of \eqn{SptAnBcongs} which we list below.
\begin{align*}
\spt( 11\cdot19^4\cdot n + 22006) &\equiv 0 \pmod{11}, \\ 
\spt( 7^4\cdot17\cdot n + 243) &\equiv 0 \pmod{17}, \\ 
\spt( 5^4\cdot19\cdot n + 99) &\equiv 0 \pmod{19}, \\ 
\spt( 13^4\cdot29\cdot n + 18583) &\equiv 0 \pmod{29}, \\ 
\spt( 29^4\cdot31\cdot n + 409532) &\equiv 0 \pmod{31}, \\ 
\spt( 5^4\cdot37\cdot n + 1349) &\equiv 0 \pmod{37}. 
\end{align*}

These explicit congruences are reminiscent of congruences for the
partition function $p(n)$ found earlier by Atkin \cite{At}, \cite{At1996},
and Weaver \cite{We}. The connection with half-integer weight Hecke eigenforms
is also analogous to what happens for the partition function. 
This was exploited by Guo and Ono
\cite{Guo-Ono}, who showed for small $\ell$ a connection between 
$\ell$-divisibility
results for the orders of certain Tate-Shafarevich groups of 
certain Tate twists of Dirichlet motives, and partition congruences.
Guo and Ono's results should extent to spt-congruences.

Recently, Folsom and Ono \cite{FolsomOno} have proved some amazing
congruences for $\spt(n)$ mod $2$ and $3$. These congruences had been observed
by the author and others. The methods of the present paper do not apply to
Folsom and Ono's results. In the present paper all congruences have 
prime modulus $\ell$ where $\ell >3$.

Recently, Bringmann \cite{Bringmann2} has shown how the generating 
function for the second rank moment gives rise to a weak Maass form of weight $3/2$.
This leads to an asymptotic formula and congruences, which in turn also implies the
result \eqn{SptAnBcongs}.

In the paper \cite{BGM} we further explore the general problem of 
congruences mod $\ell$
for rank moments, Andrews \cite{Adurf} symmetrized rank moments and full rank
functions for $k$-marked Durfee symbols for general $\ell$.
In Section \sect{conclude} we will present a 
preview of some of these results and make some concluding remarks.

\section{Congruences for $\spt(n)$ mod $5$, $7$ and $13$} \label{sec:spt5713}
The following theorem is an extension of Andrews' congruences 
\eqn{spt5}, \eqn{spt7}, and \eqn{spt13}.
\begin{theorem}
\label{thm:1}
Let $t=5$, $7$ or $13$. Then
\beq
\spt(n) \equiv -2(n + \tfrac{t^2-1}{24})\,p(n) \pmod{t},
\mylabel{eq:spt5713}
\eeq
provided $1-24n$ is not a quadratic residue mod $t$.
\end{theorem}
\begin{proof}
We proceed by considering the three cases $t=5$, $7$, $13$
separately. In each case we find the result by reducing a
known exact relation between rank and crank moments \cite{AtG}
mod $t$.

\subsection*{$t=5$}

We need \cite[(5.7),p.359]{AtG}:
\begin{align*}
N_6(n) &= \frac{2}{33}(324n^2 + 69n -10) M_2(n)
         + \frac{20}{33}(-45n + 4) M_4(n)
         + \frac{18}{11} M_6(n)  \\ 
       &\quad  + (108n^2 - 24n + 1) N_2(n).
\end{align*}
Reducing mod $5$ we obtain
\beq
N_6(n) \equiv (n^2+n) M_2(n) + 3 M_6(n) + (3n^2 + n + 1) N_2(n) 
\mylabel{eq:N6p5}
\pmod{5}.
\eeq
Since $m^{6}\equiv m^2 \pmod{5}$ we have
$$
N_{6}(n) \equiv N_2(n) \pmod{5},\qquad
M_{6}(n) \equiv M_2(n) \pmod{5}.
$$
Thus we may rewrite \eqn{N6p5}
as
\begin{align*}
2n (n+2) N_2(n) &\equiv (n + 2)(n + 4) M_2(n) \pmod{5} \\ 
&\equiv 2n (n + 2)(n + 4) p(n) \pmod{5} \qquad\mbox{(by \eqn{DyM2})},
\end{align*}
and we have
\beqs
N_2(n) \equiv (n+4) p(n) \pmod{5}\qquad\mbox{if $n\not\equiv0,3\pmod{5}$}.
\eeqs
Thus
\begin{align*}
\spt(n) &= n p(n) - \tfrac{1}{2} N_2(n)  \\ 
        &\equiv 3(n+1) p(n) \pmod{5} \qquad\mbox{if $n\not\equiv0,3\pmod{5}$},
\end{align*}
which gives Theorem \thm{1} for the case $t=5$.

\subsection*{$t=7$} This time we need \cite[(5.8),p.360]{AtG}
which is relation between the moments $M_2$, $M_4$, $M_6$, $M_8$,
$N_2$ and $N_8$. Upon reducing this relation mod $7$ and using the
fact that
$$
N_{8}(n) \equiv N_2(n) \pmod{7},\qquad
M_{8}(n) \equiv M_2(n) \pmod{7},
$$
we find that
\begin{align*}
3n(n+1)(n+5) N_2(n) &\equiv 2(n+1)(n+5)(n+6) M_2(n) \pmod{7}\\ 
&\equiv 4n(n+1)(n+5)(n+6)p(n) \pmod{7} \qquad\mbox{(by \eqn{DyM2})},
\end{align*}
and
\beqs
N_2(n) \equiv (6n+1) p(n) \pmod{7}\qquad\mbox{if $n\not\equiv0,2,6\pmod{7}$}.
\eeqs
Thus
\begin{align*}
\spt(n) &= n p(n) - \tfrac{1}{2} N_2(n)  \\ 
        &\equiv (5n+3) p(n) \pmod{7} \qquad\mbox{if $n\not\equiv0,2,6\pmod{7}$},
\end{align*}
which gives Theorem \thm{1} for the case $t=7$.

\subsection*{$t=13$} For this case we need \cite[(5.10),p.360]{AtG}
which is an exact linear
relation between the moments
$$
N_{14}, N_{12}, N_2, M_2, M_4, M_6, M_8, M_{10}, M_{12}, M_{14}.
$$
Reducing this relation mod $13$ we find that
\begin{align*}
N_{14}(n) &\equiv (4 + 4n + 12n^2 + 4n^3 + 12n^4 + 8n^6) M_2(n) \\
&\quad + (1 + 6n + 4n^2 + 2n^3 + 3n^4 + 5n^5 + 6n^6) N_2(n)
 + M_{14}(n) \pmod{13}.
\end{align*}
Using the fact that
$$
N_{14}(n) \equiv N_2(n) \pmod{13},\qquad
M_{14}(n) \equiv M_2(n) \pmod{13},\qquad
$$
we find that
\begin{align*}
&7n(n+1)(n+2)(n+5)(n+9)(n+12) N_2(n)  \\ 
&\quad\equiv 8(n+1)(n+2)(n+5)(n+9)^2(n+12) M_2(n)  \pmod{13}\\
&\quad\equiv 3n(n+1)(n+2)(n+5)(n+9)^2(n+12) p(n) \pmod{13} \qquad\mbox{(by \eqn{DyM2})},
\end{align*}
and
\beqs
N_2(n) \equiv (6n+2) p(n) \pmod{13}\qquad\mbox{if $n\not\equiv0,1,4,8,11,12
\pmod{13}$}.
\eeqs
Thus
\begin{align*}
\spt(n) &= n p(n) - \tfrac{1}{2} N_2(n)  \\ 
        &\equiv (11n+12) p(n) \pmod{13} 
\qquad\mbox{if $n\not\equiv0,1,4,8,11,12\pmod{13}$},
\end{align*}
which gives Theorem \thm{1} for the case $t=13$.
\end{proof}

\section{Modular Forms and Hecke Operators} \label{sec:modforms}

Let $N$ be a positive integer and $k$ be a nonnegative integer. We let
$M_k(N)$ (resp. $S_k(N)$) be the space of entire modular (resp. cusp)
forms of weight $k$ with respect to the modular group $\Gamma_0(N)$.
If $\chi$ is a Dirichlet character mod $N$ we let 
$M_k(N,\chi)$ (resp. $S_k(N,\chi)$) be the space of entire modular (resp. cusp)
forms of weight $k$ and character $\chi$ with respect to the modular group 
$\Gamma_0(N)$.
We define half-integral weight modular forms in the sense of Shimura \cite{Sh}.
If $N$ is divisible by $4$ 
we let $M_{k+\tfrac12}(N,\chi)$ (resp. $S_{k+\tfrac12}(N,\chi)$) be the space of 
modular (resp. cusp) forms of weight $k+\tfrac12$ and character $\chi$ with 
respect to the modular group $\Gamma_0(N)$.

Let $GL_2^{+}(\R)$ denote the group of all real $2\times2$ matrices with
positive determinant. $GL_2^{+}(\R)$ acts on the complex upper half plane 
$\mathcal{H}$ by linear fractional transformations. As in \cite{Sh}
we let $G$ denote the set of ordered pairs $\alpha,\phi(\tau)$,
where $\alpha\in GL_2^{+}(\R)$ with last row $(c\quad d)$, and $\phi$ is
a holomorphic function on $\mathcal{H}$ such that
$$
\phi^2(\tau) = s \det \alpha^{-1/2}(c\tau + d),
$$
where $\abs{s}=1$. $G$ is a group with multiplication
$$
(\alpha,\phi(\tau))(\beta,\psi(\tau))=(\alpha\beta,\phi(\beta\tau)\psi(\tau)).
$$
For a holomorphic function $f\,:\,\Hup\longrightarrow\C$ and 
$\xi=(\alpha,\phi(\tau))\in G$ we define
\beqs
f\,\mid_{k+\tfrac12}\,\xi = f\,\mid\,\xi = \phi(\tau)^{-2k-1}f(\alpha\tau).
\eeqs

For a prime $\ell$ the Hecke operator 
$T_{k,N}(\ell^2)=T(\ell^2)$ (as defined in \cite{Sh})
maps $M_{k+\tfrac12}(N,\chi)$ to itself. If $f=\sum_{n=0}^\infty a(n) q^n$
then $f\,\mid\,T(\ell^2)=\sum_{n=0}^\infty c(n) q^n$ where
\beq
c(n) = a(\ell^2n) + \chi(\ell) \leg{(-1)^kn}{\ell} \ell^{k-1} a(n)
       + \chi(\ell^2) \ell^{2k-1} a(n/\ell^2).
\mylabel{eq:Heckedef}
\eeq
We note the convention that $a(x)=0$ if $x$ is not a nonnegative integer.
Let
$$
W_N=\left(\begin{pmatrix} 0 & -1 \\ N & 0 \end{pmatrix}, N^{1/4} \sqrt{-i\tau)}\right)
\in G.
$$
Then the \textit{Fricke involution} is given by $f\mapsto f\,\mid_{k+\tfrac12}W_N$

The Dedekind eta-function is defined by
\beqs
\eta(\tau) = \exp(\pi i\tau/12) \prod_{n=1}^\infty (1 - \exp(2\pi in\tau)
= q^{1/24} \prod_{n=1}^\infty (1- q^n).
\eeqs
Then for $\tau\in\Hup$
\beq
\eta(-1/\tau) = \sqrt{-i\tau} \eta(\tau)
\mylabel{eq:etatrans}
\eeq
and $\eta(24\tau)$ is $\tfrac12$ weight cusp form in 
$S_{\tfrac12}(576,\chi_{12})$ where
\beqs
\chi_{12}(n) =
\begin{cases}
1 & \mbox{if $n\equiv\pm1\pmod{12}$,}\\
-1 & \mbox{if $n\equiv\pm5\pmod{12}$,}\\
0 & \mbox{otherwise.}                       
\end{cases}
\eeqs

\begin{prop}
\label{propo:Hecke}
Let $1\le r\le 23$ with $(r,6)=1$, $m$ be an even integer nonnegative
integer, and $\ell$ be a prime,
$\ell >3$. Define
\beqs
\mathcal{C}_{r,k} := \left\{ \eta^r(24\tau) F(24\tau)\,:\, F\in M_m(1)\right\}
\subset S_{\tfrac{r}{2}+m}(576,\chi_{12}).
\eeqs
Then the Hecke operator $T(\ell^2)$ maps $\mathcal{C}_{r,k}$ to $\mathcal{C}_{r,k}$.
\end{prop}
\begin{proof}
Suppose 
$$
g(\tau) =  \eta^r(24\tau) F(24\tau),
$$
where $F\in M_m(1)$ and $\ell$ is a prime, $\ell > 3$.
Then $g(\tau)\in S_{k+\tfrac{1}{2}}(576,\chi_{12})$ where $k=m + \tfrac{r-1}{2}$.
By \cite[(6)]{Bruinier} we have
\beqs
g \,\mid\, W_{576} \,\mid\, T(\ell^2) = 
g \,\mid\, T(\ell^2) \,\mid\, W_{576},
\eeqs
since $(\ell,576)=1$ and $\chi_{12}$ is a real character. Since $F\in M_m(1)$
we have
\beq
F(-1/\tau) = \tau^m F(\tau).
\mylabel{eq:Ftrans}
\eeq
Hence using \eqn{etatrans} and \eqn{Ftrans} we find that
\begin{align*}
g\,\mid\, W_{576} 
& = \left(\sqrt{24} \sqrt{-i\tau}\right)^{-2k-1} g(-1/(576\tau)) \\
& = \left(\sqrt{24} \sqrt{-i\tau}\right)^{-2k-1} 
\eta^r(-1/(24\tau)) F(-1/(24\tau)\\
&= i^m g. 
\end{align*}
Now let
\beqs
H(\tau) = \frac{ g\,\mid\,T(\ell^2)}{\eta^r(24\tau)}.
\eeqs
Then
\begin{align}
H(-1/(576\tau))
&= \left(\sqrt{24} \sqrt{-i\tau}\right)^{2k+1}
   \frac{g\,\mid\,T(\ell^2)\,\mid\,W_{576}}
        {\eta^2(-1/(24\tau))} \mylabel{eq:Htrans}\\
&=\left(\sqrt{24} \sqrt{-i\tau}\right)^{2k+1}
   \frac{g\,\mid\,W_{576}\,\mid\,T(\ell^2)}
        {(-24i\tau)^{r/2}\eta^2(-1/(24\tau))} \nonumber\\
&=i^m (-24i\tau)^{k+\tfrac{1-r}{2}}\frac{ g\,\mid\,T(\ell^2)}{\eta^r(24\tau)}
\nonumber\\
&=(24\tau)^m H(\tau).
\nonumber
\end{align}
We note the exponents of $q$ in the $q$-expansion of $g(\tau)$
are congruent to $r\pmod{24}$. Since $\ell^2\equiv1\pmod{24}$
we see by \eqn{Heckedef} that the exponents of $q$ in the $q$-expansion of 
$g(\tau)\,\mid\,T(\ell^2)$ are also congruent to $r\pmod{24}$. Since
$$
\eta^r(24\tau) = q^r + \cdots
$$
and $\eta(\tau)$ is nonzero in $\Hup$
we see that 
$$
H(\tau) = K(24\tau),
$$
for some function $K(\tau)$ holomorphic function on $\Hup$ that 
satisfies
$$
K(\tau + 1) = K(\tau).
$$
From \eqn{Htrans} we have
$$
K(-1/\tau) = H(-1/(24\tau)) = \tau^m H(\tau/24) = \tau^m K(\tau).
$$
Hence $K(\tau)\in M_m(1)$ and
$$
g\,\mid\,T(\ell^2) = \eta^r(24\tau) K(24\tau) \in \mathcal{C}_{r,k}.
$$
\end{proof}
\begin{cor}
\mylabel{cor:Heckecor}
Let $1\le r\le 23$ with $(r,6)=1$, and $m$ be an even integer nonnegative
integer. If $\dim M_m(1)=1$ and $0\ne F(\tau)\in M_m(1)$, then the function
$$
g(\tau) = \eta^r(24\tau) F(24\tau)
$$
is a Hecke eigenform for 
$S_{\tfrac{r}{2}+m}(576,\chi_{12})$.
\end{cor}

We define the slash operator for modular forms of
integer weight. Let $k\in\Z$.
For a holomorphic function $f\,:\,\Hup\longrightarrow\C$ 
and $\alpha\begin{pmatrix} a & b \\ c & d \end{pmatrix}\in GL_2^{+}(\R)$ 
we define
\beqs
f\,\mid_{k}\,\alpha = f\,\mid\,\alpha = 
(\det \alpha)^{\tfrac{k}{2}} (c\tau + d)^{-k} f(\alpha \tau).
\eeqs
For $m$ a positive integer we define the operator $U(m)$ by
its action on formal power series
\beqs
\left( \sum_{n=0}^\infty a(n) q^n\right) \,\mid\,U(m)
= \sum_{n=0}^\infty a(mn) q^n.
\eeqs

We will also need the Hecke operator $T(\ell)$ on $M_k(1)$.
For a positive integer $m$ the Hecke operator $T(m)$ (as defined
in \cite{Ap}) maps $M_k(1)$ and $S_k(1)$ to themselves. 
Let $\ell$ be prime.
If $f=\sum_{n=0}^\infty a(n) q^n$
then $f\,\mid\,T(\ell)=\sum_{n=0}^\infty c(n) q^n$ where
\beqs
c(n) = a(\ell n) 
       + \ell^{k-1} a(n/\ell).
\eeqs
Hence, if the coefficients $a(n)$ are integers we have
\beqs
f\mid T(\ell) \equiv f \mid U(\ell) \pmod{\ell^{k-1}}.
\eeqs

Following Chua \cite{Chua} we  define
\beq
h_{\ell}(\tau) = \left(\eta(\tau) \eta(\ell\tau)\right)^{\ell-1},
\mylabel{eq:hell}
\eeq
for $\ell>3$ prime. The following proposition is basically an extension of
a result of Chua \cite{Chua}. The $F=1$ case follows from 
Lemmas 2.1 and 2.2 in \cite{Chua}.

\begin{prop}
\mylabel{propo:Chuaext}
Let $\ell>3$ be prime, $k$ be a nonnegative even integer.
If $F(\tau)\in M_k(1)$ then
\beq
h_{\ell}(\tau) F(\tau) \mid U(\ell) + 
(-1)^{(\ell-1)/2} \ell^{k + (\ell-1)/2 -1} h_{\ell}(\tau) F(\ell \tau)
\in S_{k+\ell-1}(1).
\mylabel{eq:helltrans}
\eeq
\end{prop}
\begin{proof}
Suppose $\ell>3$ is prime, $k$ is an nonnegative even integer, and
$F(\tau)\in M_k(1)$.          
Chua \cite{Chua} showed that $h_{\ell}(\tau)\in M_{\ell-1}(\ell)$.
Actually $h_{\ell}(\tau)\in S_{\ell-1}(\ell)$ since clearly
$h_{\ell}$ is zero at the cusp $i\infty$ and zero at the cusp $0$
by \eqn{helltrans2} below.
Hence 
\beqs
g_{\ell}(\tau) = h_{\ell}(\tau) F(\tau) \in S_{k+\ell-1}(\ell).
\eeqs
By \cite[Lemma 17(iii), p.144]{AtkinLehner},
\beq
g_{\ell}(\tau) \mid U(\ell) + 
\ell^{k+\ell-1} g_{\ell}(\tau) \mid W_{\ell}\in S_{k+\ell-1},
\mylabel{eq:gtrans}
\eeq
where
\beqs
W_{\ell}=\begin{pmatrix} 0 & -1 \\ \ell & 0\end{pmatrix}.
\eeqs
By \eqn{etatrans} we have
\beq
h_{\ell}\left(-\tfrac{1}{\ell\tau}\right)
 = (-1)^{(\ell-1)/2} \ell^{(\ell-1)/2} \tau^{\ell-1} h_{\ell}(\tau).
\mylabel{eq:helltrans2}
\eeq
Since $F(\tau)\in M_k(1)$, we have
\beq
F\left(-\tfrac{1}{\ell\tau}\right) = \ell^k \tau^k F(\ell\tau).
\mylabel{eq:Ftrans2}
\eeq
Now by \eqn{helltrans2} and \eqn{Ftrans2} we have
\beq
g_{\ell} \mid_{k+\ell-1} W_{\ell}
= \ell^{-(k+\ell-1)/2} \tau^{-k - \ell + 1}
 g_{\ell}\left(-\tfrac{1}{\ell\tau}\right)
 = (-1)^{(\ell-1)/2} \ell^{k/2} h_{\ell}(\tau) F(\ell\tau).
\mylabel{eq:gtrans2}
\eeq
The result \eqn{helltrans} follows from \eqn{gtrans} and \eqn{gtrans2}. 
\end{proof}

Let $F\in M_k(1)$ and suppose the $q$-expansion of $F(\tau)$
has integer coefficients. We call such a form an \textit{integral}
modular form.
We define $p(F,n)$ by
\beqs
\sum_{n=0}^\infty p(F,n) q^n = \frac{F(\tau)}{\prod_{n=1}^\infty (1-q^n)}.
\eeqs

We prove a straightforward generalization of a
theorem due to Chua \cite[Theorem 1.1]{Chua}.
We note that Chua's result was extended to prime power moduli
by Ahlgren and Boylan \cite{AB}.

\begin{theorem}
\mylabel{thm:pFcong} 
Suppose $\ell>3$ is prime, $k$ is an nonnegative even integer, and
$F(\tau)\in M_k(1)$ is an integral modular form. Define $1\le \beta_{\ell}\le \ell-1$
such that $24\beta_{\ell}\equiv1\pmod{\ell}$ and let
\beqs
r_{\ell} = \frac{24\beta_{\ell}-1}{\ell},
\quad
\lambda_{\ell} = \frac{\ell^2 + 24\beta_{\ell}-1}{24\ell}.
\eeqs
Then
\beqs
\sum_{n=0}^\infty p(F,\ell n +\beta_{\ell}) q^{24n+r_{\ell}}
\equiv \eta^{r_{\ell}}(24\tau) G_{\ell,F}(24\tau) \pmod{\ell},
\eeqs
for some integral
modular form $G_{\ell,F}(\tau)\in M_{k+\ell-1 - 12\lambda_{\ell}}(1)$.
\end{theorem}
\begin{proof}
Suppose $\ell>3$ is prime, $k$ is an nonnegative even integer, and
$F(\tau)\in M_k(1)$. 
\begin{align*}
\sum_{n=0}^\infty p(F,n) q^{24n-1} & = \frac{F(24\tau)}{\eta(24\tau)}
 \\
&\equiv  \frac{F(24\tau)}{\eta(24\tau)} \, 
 \frac{\eta^{\ell}(24\tau)}{\eta(24\ell\tau)} \pmod{\ell}
\\
&\equiv \frac{h_{\ell}(24\tau) F(24\tau)}{\eta^{\ell}(24\ell\tau)}
\pmod{\ell},
\end{align*}
where $h_{\ell}(\tau)$ is defined in \eqn{hell}.
Applying the $U(\ell)$ operator to both sides we obtain
\beqs
\sum_{n=0}^\infty p(F,\ell n+\beta_{\ell}) q^{24n+r_\ell} 
\equiv \frac{h_{\ell}(24\tau) F(24\tau)\mid U(\ell)}{\eta^{\ell}(24\tau)}
\pmod{\ell}.
\eeqs
By Proposition \propo{Chuaext} there is cusp form $j_{\ell}(\tau)
\in M_{k+\ell -1}$ such that
\begin{align}
j_{\ell}(\tau) &=
h_{\ell}(\tau) F(\tau) \mid U(\ell) + 
(-1)^{(\ell-1)/2} \ell^{k + (\ell-1)/2 -1} h_{\ell}(\tau) F(\ell \tau)
\mylabel{eq:jell}\\
&=\left( c_1 q^{ (\ell^2-1)/24} + \cdots \right)\mid U(\ell) 
+ O( q^{ (\ell^2-1)/24}),
\nonumber\\
&= c_2 q^{\lambda_{\ell}} + \cdots
\nonumber
\end{align}
for some constants $c_1$, $c_2$ since for $\ell\ge5$,
$\lambda_{\ell} \le \frac{\ell^2-1}{24}$. It follows that
\beq
j_{\ell}(\tau) = \Delta^{\lambda_\ell}(\tau) G_{\ell,F}(\tau),
\mylabel{eq:jellid}
\eeq
for some $G_{\ell,F}(\tau)\in M_{k +\ell -1 - 12\lambda_\ell}(1)$.
Here as usual
\beqs
\Delta(\tau) = \eta^{24}(\tau).
\eeqs
Since $\ell\ge5$, $k + (\ell-1)/2 -1\ge 1$, and
\beqs
h_{\ell}(\tau) F(\tau) \mid U(\ell) 
\equiv  \Delta^{\lambda_\ell}(\tau) G_{\ell,F}(\tau)\pmod{\ell},
\eeqs
by \eqn{jell} and \eqn{jellid}. Since $(\ell,24)=1$ we have
\beqs
\sum_{n=0}^\infty p(F,\ell n+\beta_{\ell}) q^{24n+r_\ell}
\equiv \frac{ \Delta^{\lambda_\ell}(24\tau) G_{\ell,F}(24\tau)}{\eta^{\ell}(24\tau)}
\equiv 
 \eta^{r_{\ell}}(24\tau) G_{\ell,F}(24\tau) \pmod{\ell},
\eeqs
since $24\lambda_\ell -\ell = r_\ell$.
\end{proof}

\section{Congruences for {\rm $\spt(n)$} mod $11$} \label{sec:spt11}

We need some results for the crank mod $11$. For $t\ge2$ and
$0\le r \le t-1$ let $M(r,t,n)$ denote the number of partitions
of $n$ with crank congruent to $r$ mod $t$.
Then for $t=5$, $t=7$, or $t=11$
\beqs
M(r,t,n) = \frac{1}{t} \, p(n), \qquad 0\le r \le t-1;
\eeqs
for all $n$ satisfying $24n\equiv 1 \pmod{t}$.
See \cite{AG}, \cite{G88a}.
These combinatorial results immediately imply Ramanujan's
partition congruences \eqn{p5}, \eqn{p7}, \eqn{p11}.
Here we need the $t=11$ case
\beq
M(r,11,11n+6) = \frac{1}{11} \, p(11n+6), \qquad 0\le r \le 10.
\mylabel{eq:Mrels11}
\eeq
Let
\begin{equation*}
E(q) = \prod_{n=1}^\infty (1-q^n).
\end{equation*}
For $r\ge1$ define $p_r(n)$ by
\beqs
\sum_{n=0}^\infty p_r(n) q^n = E(q)^r = \prod_{n=1}^\infty (1-q^n)^r,
\eeqs
so that
\beqs
\sum_{n=1}^\infty p_{23}(n-1) q^n = qE(q)^{23} = q\prod_{n=1}^\infty (1-q^n)^{23}.
\eeqs
We need \cite[(5.24),pp.362-3]{AtG} which is an exact linear relation
between rank-crank moments and the function $p_{23}(n-1)$.
Upon reducing this relation mod $11$ and            
using the fact that
$$
N_{12}(n) \equiv N_2(n) \pmod{11},\qquad
M_{12}(n) \equiv M_2(n) \pmod{11},\qquad
$$
we find that
\begin{align*}
p_{23}(11n+5)
&\equiv 4 N_2(11n+6) + 2 M_2(11n+6) \\ 
&\quad     + M_4(11n+6) + M_{6}(11n+6) + 10 M_8(11n+6) 
\pmod{11}.
\end{align*}
By \eqn{Mrels11} we have
\begin{align*}
M_k(11n+6) &\equiv \sum_{m=1}^{10} m^k M(m,11,11n+6)\pmod{11}  \\ 
&\equiv M(1,11,11n+6) \sum_{m=1}^{10} m^k \pmod{11}\\
&\equiv 0\pmod{11}, \qquad\mbox{if $k=2,4,6,8$.}
\end{align*}
Hence
\beq
N_2(11n+6) \equiv 3 p_{23}(11n+5) \pmod{11}
\mylabel{eq:N2m11}
\eeq
and
\begin{align*}
\spt(11n+6) &= (11n+6) p(11n+6) - \tfrac{1}{2}N_2(11n+6)   \\ 
&\equiv 4 p_{23}(11n+5) \pmod{11}.
\end{align*}
Now by the Children's Binomial Theorem we have
\begin{align*}
\prod_{n=1}(1-q^n)^{23} &= \prod_{n=1}(1-q^n)^{22} \prod_{n=1}(1-q^n)
\\
&\equiv \prod_{n=1}(1-q^{11n})^2 \prod_{n=1}(1-q^n) \pmod{11}.
\end{align*}
We need Euler's pentagonal number theorem
\beqs
\prod_{n=1}^\infty (1 - q^n) = \sum_{n=-\infty}^\infty (-1)^n q^{n(3n-1)/2}.
\eeqs
Since $n(3n-1)/2\equiv5\pmod{11}$ if and only if $n\equiv2\pmod{11}$
we find that
\begin{align}
\sum_{n=0}^\infty p_{23}(11n+5) q^n 
&\equiv \prod_{n=1}(1-q^n)^2 \prod_{n=1}(1-q^{11n})\pmod{11}\mylabel{eq:p23m11b}\\
&\equiv \prod_{n=1}(1-q^n)^{13} \pmod{11}.
\nonumber
\end{align}
Hence
\beqs
\spt(11n+6) \equiv 4 p_{13}(n) \pmod{11}.
\eeqs
Now let
\beqs
F_{13}(\tau) = \eta^{13}(24\tau) =
\sum_{n\ge13} b_{13}(n) q^n 
\in S_{\tfrac{13}{2}}(576,\chi_{12}),
\eeqs
where
\beqs
b_{13}(n)=p_{13}\left(\frac{n-13}{24}\right).
\eeqs
Then by Theorem \propo{Hecke}, $F_{13}$ is a Hecke eigenform so
\beq
F_{13} \,\mid\, T(\ell^2) = \lambda_{\ell}\, F_{13}, \qquad(\ell>3),
\mylabel{eq:FT13}
\eeq
where $\lambda_{\ell}\in\Z$. 
We look for the smallest eigenvalue which is a multiple of $11$.
After some calculation we find that
\beqs
\lambda_{19} = - 2901404 = - 2^2\cdot 11 \cdot 23 \cdot 47 \cdot 61
\equiv 0 \pmod{11},
\eeqs
By \eqn{FT13} and \eqn{Heckedef} we have
\beqs
b_{13}(19^3m) = \lambda_{19} \,b_{13}(19m) \equiv 0 \pmod{11}
\qquad\mbox{when $(m,19)=1$.}
\eeqs
We want
$$
19^3m\equiv 13 \pmod{24},\quad\mbox{and}\quad (m,19)=1,
$$
and take
$$
m = 24\cdot 19 k + 7,\quad n = 19^3(24\cdot 19 k + 7)
$$
so that
$$
\frac{n-13}{24} = 19^4k + 2000.
$$
Hence
$$
p_{13}(19^4n + 2000) \equiv 0 \pmod{11},
$$
and
$$
\spt(19^4\cdot 11n + 22006) \equiv 0 \pmod{11}.
$$

\section{Explicit Congruences for the rank mod $11$} 
\label{sec:explicitrankcongs}

In this section we find explicit congruences for the rank mod $11$.
From \eqn{N2m11} and \eqn{p23m11b} we find that the generating function
for $N_2(11n+6)$ is congruent to $3 E(q)^{13}\pmod{11}$.
In the following theorem we give congruences mod $11$ for
other rank moments $N_{2k}(11n+6)$ in terms of half-integer
weight cusp forms. 

Following \cite{AtG} and Ramanujan \cite[p.163]{Ram} we define
\beqs
\Phi_j = \Phi_j(q) 
= \sum_{n=1}^\infty \frac{n^j q^n}{1-q^n} = \sum_{m,n\ge1} n^j q^{nm}
= \sum_{n=1}^\infty \sigma_j(n) q^n,
\eeqs
for $j\ge1$ odd and where $\sigma_j(n) = \sum_{d\mid n} d^j$.
For $n\ge2$ even we define the Eisenstein series
\beqs
E_n(\tau) 
= 1 - \frac{2n}{B_n} \Phi_{n-1}(q),
\eeqs
where $q=\exp(2\pi i\tau)$, $\Im \tau >0$ and $B_n$ is the $n$-th
Bernoulli number. We note that $E_n(\tau)\in M_n(1)$ for $n\ge 4$.

\begin{theorem}
\mylabel{thm:rankmomsmod11}
\begin{align}
\sum_{n=0}^{\infty }{N}_{{2}} \left( 11 n+6 \right) {q}^{n}
&\equiv 3 {E}^{13}(q)
\pmod{11},
\mylabel{eq:N2mod11}\\
\sum_{n=0}^{\infty }{N}_{{4}} \left( 11 n+6 \right) {q}^{n}
&\equiv 7 {E}^{13}(q)
\pmod{11},
\mylabel{eq:N4mod11}\\
\sum_{n=0}^{\infty }{N}_{{6}} \left( 11 n+6 \right) {q}^{n}
&\equiv {E}^{13}(q) \left( 4+E_{{4}}(\tau) \right) 
\pmod{11}, 
\mylabel{eq:N6mod11}\\
\sum_{n=0}^{\infty }{N}_{{8}} \left( 11 n+6 \right) {q}^{n}
&\equiv{E}^{13}(q) \left( 5+6 E_{{4}}(\tau)+6 E_{{6}}(\tau) \right)
\pmod{11}, 
\mylabel{eq:N8mod11}\\
\sum_{n=0}^{\infty }{N}_{{10}} \left( 11 n+6 \right) {q}^{n}
&\equiv{E}^{13}(q) 
\left( 5+4 E_{{4}}(\tau)+6 E_{{6}}(\tau)+6 {E_{{4}}}^{2}(\tau) \right) 
\pmod{11}.
\mylabel{eq:N10mod11}
\end{align}
\end{theorem}

It is clear the each rank moment $N_{2k}(11n+6)\pmod{11}$ can
be written in terms the rank functions $N(r,11,11n+6)$, ($0 \le r \le 5$).
These relations may be inverted to find each $N(r,11,11n+6)$
in terms of rank moments mod $11$.
We note that only $r\le 5$ is needed since $N(-m,n)=N(m,n)$ and
$N(11-r,11,n)=N(r,11,n)$.

\begin{cor}
\mylabel{cor:rankcongsmod11}
\begin{align*}
\sum_{n=0}^{\infty }{N} \left( 0,11,11 n+6 \right) {q}^{n}
&\equiv {E}^{13}(q) 
\left( 6+7 {E_{{4}}}(\tau)+5 {E_{{6}}}(\tau)+5 {E_{{4}}}^{2}(\tau) \right) 
\pmod{11}, \\
\sum_{n=0}^{\infty }{N} \left( 1,11,11 n+6 \right) {q}^{n}
&\equiv {E}^{13}(q) 
\left( 9+10 {E_{{6}}}(\tau)+5 {E_{{4}}}^{2}(\tau) \right) 
\pmod{11}, \\
\sum_{n=0}^{\infty }{N} \left( 2,11,11 n+6 \right) {q}^{n}
&\equiv {E}^{13}(q) 
\left( 4+3 {E_{{6}}}(\tau)+5 {E_{{4}}}^{2}(\tau) \right) 
\pmod{11}, \\
\sum_{n=0}^{\infty }{N} \left( 3,11,11 n+6 \right) {q}^{n} 
&\equiv {E}^{13}(q) 
\left( 8+4 {E_{{4}}}(\tau)+6 {E_{{6}}}(\tau)+5 {E_{{4}}}^{2}(\tau) \right) 
\pmod{11}, \\
\sum_{n=0}^{\infty }{N} \left( 4,11,11 n+6 \right) {q}^{n}
&\equiv {E}^{13}(q) 
\left( 2+7 {E_{{4}}}(\tau)+8 {E_{{6}}}(\tau)+5 {E_{{4}}}^{2}(\tau) \right) 
\pmod{11}, \\
\sum_{n=0}^{\infty }{N} \left( 5,11,11 n+6 \right) {q}^{n}
&\equiv {E}^{13}(q) 
\left( 7+2 {E_{{4}}}(\tau)+9 {E_{{6}}}(\tau)+5 {E_{{4}}}^{2}(\tau) \right) 
\pmod{11}.
\end{align*}
\end{cor}

Although the half-integer modular forms appearing
in the previous theorem have different weights each one is
a Hecke eigenform (in its corresponding space) in view
of Corollary \corol{Heckecor}. As a result we find following
rank congruences. 

\begin{cor}
\mylabel{cor:exprankcongsmod11}
\begin{align}
N(r,11,5^4\cdot11\cdot19^4\cdot n + 4322599) &\equiv 0 \pmod{11},
\mylabel{eq:nicecong11a}\\
N(r,11,11^2\cdot19^4\cdot n + 172904) &\equiv 0 \pmod{11},
\mylabel{eq:nicecong11b}
\end{align}
for all $0\le r \le 10$.
\end{cor}

Before we can prove Theorem \thm{rankmomsmod11} we need to
recall some results on rank and crank moments from \cite{AtG}.
We define the  moment generating functions.
\begin{align*}
R_k(q) &= \sum_{n\ge0} N_k(n) q^n,\\
C_k(q) &= \sum_{n\ge0} M_k(n) q^n, 
\end{align*}
for $k$ even. 
Define
\beqs
P = P(q) = \prod_{n=1}^\infty \frac{1}{(1-q^n)},
\eeqs
for $\abs{q}<1$. In \cite{AtG} we proved the following recurrence
for crank moments
\beqs
C_{2n} = 2\sum_{j=1}^{n - 1}
      \binom{2n-1}{2j-1}\,\Phi_{2j-1}\,C_{2n-2j}
      + 2 \Phi_{2n-1}\,P,
\eeqs
which implies that there are integers 
$\alpha_{a_1,a_2,\dots,\alpha_n}$ such that
\beq
C_{2n} = 
2\, P\, \sum_{a_1+2a_2 + \cdots + n a_n = n} \alpha_{a_1,a_2,\dots,\alpha_n} \Phi_1^{a_1} \Phi_3^{a_2}
\cdots \Phi_{2n-1}^{a_n},
\mylabel{eq:C2nform}
\eeq
In \cite{AtG} we obtained the following identity 
\begin{align}
&\sum_{i=0}^{k-1}\binom{2k}{2i}
\sum_{\substack{\alpha+\beta+\gamma=2k-2i\\ 
      \mbox{$\alpha$, $\beta$, $\gamma$ even $\ge0$}}}
\binom{2k-2i}{\alpha,\beta,\gamma}\,C_\alpha\,C_\beta\,C_\gamma\,P^{-2}
- 3\left(2^{2k-1}-1\right) C_2 \mylabel{eq:rcrel}\\
&= \frac{1}{2}(2k-1)(2k-2) R_{2k}
+ 6\sum_{i=1}^{k-1} \binom{2k}{2i}\left(2^{2i-1}-1\right) 
\delta_q(R_{2k-2i})\nonumber
\\
&\quad + \sum_{i=1}^{k-1}\left[
\binom{2k}{2i+2}\left(2^{2i+1}-1\right) - 2^{2i}\binom{2k}{2i+1} 
+ \binom{2k}{2i}
\right] R_{2k-2i}.
\nonumber
\end{align}
Here $\delta_q$ is the differential operator
$$
\delta_q = q\, \frac{d}{dq}.
$$

\bigskip

\noindent
\textit{Proof of Theorem \thm{rankmomsmod11}}. 
As noted above, the first congruence
\eqn{N2mod11} follows from \eqn{N2m11} and \eqn{p23m11b}.
To attack \eqn{N4mod11}--\eqn{N10mod11}, we first use 
\eqn{rcrel}, and \eqn{C2nform} rewriting each $\Phi_{2j-1}$
in terms of $E_2(\tau)$, $E_4(\tau)$ and $E_6(\tau)$. 
We define the operator
$U_{11}^{*}$
which acts on $q$-series by
\beqs
U_{11}^{*}\left(\sum_{n=0}^\infty a(n) q^n \right)
=
\sum_{n=0}^\infty a(11n + 6) q^n.
\eeqs
We find that
\begin{align}
U_{11}^{*}\left({R_4}\right) &\equiv
U_{11}^{*}\left( 6 R_2 + P H^{(4)}\right) \pmod{11},
\mylabel{eq:U11R4}\\
U_{11}^{*}\left({R_6}\right) &\equiv
U_{11}^{*}\left( 5 R_2 + P H^{(6)}\right) \pmod{11},
\mylabel{eq:U11R6}\\
U_{11}^{*}\left({R_8}\right) &\equiv
U_{11}^{*}\left( 9 R_2 + P H^{(8)}\right) \pmod{11},
\mylabel{eq:U11R8}\\
U_{11}^{*}\left({R_{10}}\right) &\equiv
U_{11}^{*}\left( 9 R_2 + \Psi + P H^{(10)}\right) \pmod{11},
\mylabel{eq:U11R10}
\end{align}
where
\begin{align*}
H^{(4)}
&=
9 + 5\,E_{{2}} + 10\,E_{{4}} + 9\,{E_{{2}}}^{2} \\
H^{(6)}
&=
2 + 7\,E_{{2}} + 7\,E_{{4}} + 3\,{E_{{2}}}^{2} + 8\,E_{{6}} + 8\,E_{{4}}E_{{2}}
 + 9\,{E_{{2}}}^{3} \\
H^{(8)}
&=
6\,E_{{2}} + {E_{{2}}}^{2} + 6\,E_{{4}} + 10\,{E_{{2}}}^{3} + 4\,E_{{6}}
 + 4\,E_{{4}}E_{{2}} + 9\,{E_{{2}}}^{4} + 5\,{E_{{2}}}^{2}E_{{4}} \\
&\qquad + 10\,E_{{6}}E_{{2}}
\\
H^{(10)}
&=
8 + 4\,E_{{2}} + 8\,{E_{{2}}}^{2} + 4\,E_{{4}} + 3\,{E_{{2}}}^{3}
 + 10\,E_{{6}} + 10\,E_{{4}}E_{{2}} + 9\,{E_{{2}}}^{4} + 5\,{E_{{2}}}^{2}E_{{4}}
\\
&\qquad + 10\,E_{{6}}E_{{2}} + 4\,{E_{{2}}}^{5} + 10\,E_{{4}}E_{{6}}
 + 5\,{E_{{2}}}^{2}E_{{6}} + 9\,{E_{{2}}}^{3}E_{{4}}
\end{align*}
and
\beq
\Psi = \tfrac{1}{11} P (1 - E_4 E_6).
\mylabel{eq:Psidef}
\eeq
We call a function $11$-integral if the coefficients in its
$q$-expansion are rational numbers with bounded denominators
relatively prime to $11$.
We note that all functions involved are $11$-integral so that
each congruence is well-defined.
For example, $\Psi$ is $11$-integral since
\begin{equation}
E_{4} E_{6} = E_{10} = 1 - 264\,\Phi_9 \equiv 1 \pmod{11}.
\mylabel{eq:Psicong}
\end{equation}

For each $m$ we write
\beqs
H^{(m)} = H^{(m)}_0 +  H^{(m)}_2+  H^{(m)}_4 + H^{(m)}_6 +  H^{(m)}_8,
\eeqs
where each $H^{(m)}_j$ is congruent mod $11$ to an $11$-integral modular 
form in $M_{\kappa(m,j)}(1)$. Here the weight 
$\kappa=\kappa(m,j)\equiv j\pmod{10}$.  Here we have used
\cite[Theorem 2,p.22]{SWD}. We note that
\beqs
E_2 \equiv E_{12} \pmod{11},\qquad\mbox{and}\qquad
E_{10} \equiv 1 \pmod{11}.
\eeqs
By Theorem \thm{pFcong}, each $U_{11}^{*}(P H^{(m)}_j)$ is congruent
to a function $E^{13}(q) G_{m,j}(\tau)$ where $G_{m,j}(\tau)$
to a modular form of weight
$\kappa(m,j)-2$.
If $H^{(m)}_j$ is nonzero we 
prove a   congruence for $U_{11}^{*}(P H^{(m)}_j)$
by checking enough coefficients of $q^n$
(i.e. $n\le \lfloor\tfrac{\kappa(m,j)-2}{12}\rfloor \le 5$).
This is a standard argument.  See for example \cite{SWD}, \cite{St}.
For example, for $m=10$ we have
\beqs
H^{(10)}=
H^{(10)}_{0}
+H^{(10)}_{2}
+H^{(10)}_{4}
+H^{(10)}_{6}
+H^{(10)}_{8}
\eeqs
where
\begin{align*}
H^{(10)}_{0} &= 8
 + 4\,{E_{{2}}}^{5}
 + 10\,E_{{4}}E_{{6}}
 + 5\,{E_{{2}}}^{2}E_{{6}}
 + 9\,{E_{{2}}}^{3}E_{{4}}
\qquad\mbox{($\kappa=60$)}  \\
H^{(10)}_{2} &= 4\,E_{{2}}
\qquad\mbox{($\kappa=12$)}  \\
H^{(10)}_{4} &= 8\,{E_{{2}}}^{2}
 + 4\,E_{{4}}
\qquad\mbox{($\kappa=24$)}  \\
H^{(10)}_{6} &= 3\,{E_{{2}}}^{3}
 + 10\,E_{{6}}
 + 10\,E_{{4}}E_{{2}}
\qquad\mbox{($\kappa=36$)}  \\
H^{(10)}_{8} &= 9\,{E_{{2}}}^{4}
 + 5\,{E_{{2}}}^{2}E_{{4}}
 + 10\,E_{{6}}E_{{2}}
\qquad\mbox{($\kappa=48$)}
\end{align*}
We prove the following congruences 
\begin{align*}
U_{11}^{*}(P H^{(10)}_{0}) &\equiv 4\, E^{13}(q){E_{{4}}}^{2}
\pmod{11},\\
U_{11}^{*}(P H^{(10)}_{2}) &\equiv 0
\pmod{11},\\
U_{11}^{*}(P H^{(10)}_{4}) &\equiv 0
\pmod{11},\\
U_{11}^{*}(P H^{(10)}_{6}) &\equiv 4\, E^{13}(q)E_{{4}}
\pmod{11},\\
U_{11}^{*}(P H^{(10)}_{8}) &\equiv 6\, E^{13}(q)E_{{6}}
\pmod{11},
\end{align*}
by checking enough terms.
Hence    
\beq
U_{11}^{*}(P H^{(10)} ) \equiv  E^{13}(q) \left( 4\,{E_{{4}}}^{2}+
4\,E_{{4}}+6\,E_{{6}}\right) \pmod{11}.
\mylabel{eq:UPH10}
\eeq
Similarly, we find that
\begin{align}
U_{11}^{*}(P H^{(4)} ) &\equiv  0 \pmod{11}, \mylabel{eq:UPH4}\\
U_{11}^{*}(P H^{(6)} ) &\equiv  E^{13}(q) \, E_4 \pmod{11},
\mylabel{eq:UPH6}\\
U_{11}^{*}(P H^{(8)} ) &\equiv  E^{13}(q) \left( 6\,E_4 
+6 \, E_6\right) \pmod{11}.
\mylabel{eq:UPH8}
\end{align}

To prove \eqn{N10mod11} we also need to compute $U_{11}^{*}(\Psi)$ mod $11$.
Since 
\beqs
\left(\frac{E(q)^{11}}{E(q^{11})}\right)^{11}
\equiv 1 \pmod{11^2},
\eeqs
we have
\beq
q U_{11}^{*}(P(1 - E_4 E_6) ) \equiv \frac{ \Delta^5(1 - E_4 E_6)\mid U_{11} }
                                          {E(q)^{11}}
\pmod{11^2}.
\mylabel{eq:U11Psi}
\eeq
Now $\Delta^5\in S_{60}(1)$ and $E_4 E_6 \Delta^5\in S_{70}(1)$.
The set
\beqs
\mathcal{B}_1 = \{\Delta{E_{{4}}}^{12},{\Delta}^{2}{E_{{6}}}^{6},
{\Delta}^{3}{E_{{6}}}^{4},{\Delta}^{4}{E_{{6}}}^{2},{\Delta}^{5}\}
\eeqs
is a basis for $S_{60}(1)$ and multiplying each element by
$E_4 E_6$ gives a basis $\mathcal{B}_2$ for $S_{70}(1)$. We have
\begin{align*}
&\Delta^5 \mid U_{11} \equiv \Delta^5 \mid_{60} T(11) \pmod{11^2},\\
&E_4 E_6 \Delta^5 \mid U_{11} \equiv \Delta^5 \mid_{72} T(11) \pmod{11^2}.
\end{align*}
We express $\Delta^5 \mid_{60} T(11)$ in terms of $\mathcal{B}_1$,
$E_4 E_6 \Delta^5 \mid_{72} T(11)$ in terms of $\mathcal{B}_2$, and reduce 
mod $11^2$ to find
\begin{align*}
\Delta^5(1 - E_4 E_6)\mid U_{11}
&\equiv 11 
\Delta\, \left( 1+E_{{4}}E_{{6}} \right)  \left( {E_{{4}}}^{12}+10\,{\Delta}^{4}+
6\,{\Delta}^{2}{E_{{6}}}^{4}+7\,{\Delta}^{3}{E_{{6}}}^{2}
+9\,\Delta\,{E_{{6}}}^{6} \right)\\
&\equiv 22 E_4^2 \Delta \pmod{11^2},
\end{align*}
by using 
\beqs
\Delta = \tfrac{1}{1728}(E_4^3 - E_6^2),
\eeqs
and \eqn{Psicong}
It follows from \eqn{U11Psi} that
\beq
U_{11}^{*}(\Psi) \equiv 2 E^{13}(q) E_4^2 \pmod{11}.
\mylabel{eq:U11Psi2}
\eeq
Finally, \eqn{N10mod11} follows from \eqn{N2mod11}, \eqn{U11R10},
\eqn{UPH10} and \eqn{U11Psi2}. This completes the
proof of Theorem \thm{rankmomsmod11}.

\bigskip

\noindent
\textit{Proof of Corollary \corol{exprankcongsmod11}}. 
From Corollary \corol{rankcongsmod11} we see that for each $0\le r\le 10$,
that there are integers $a_r$, $b_r$, $c_r$, and $d_r$ such
\beq
\sum_{n=0}^\infty
N(r,11,\tfrac{1}{24}(11 n + 1)) q^n
\equiv \eta^{13}(24\tau)
(a_r + b_r E_4(24\tau) + c_r E_6(24\tau) + d_r E_4^2(24\tau))
\pmod{11}.
\mylabel{eq:Nrmod11}
\eeq
By Corollary \corol{Heckecor} we note that each of the functions
\beqs
\eta^{13}(24\tau),\quad
\eta^{13}(24\tau) E_4(24\tau),\quad
\eta^{13}(24\tau) E_6(24\tau),\quad
\eta^{13}(24\tau) E_4^2(24\tau),
\eeqs
is a Hecke eigenform in its corresponding space;
i.e. $S_{w}(576,\chi_{12})$ for $w=\tfrac{13}{2}$,
$w=\tfrac{13}{2}$,
$w=\tfrac{21}{2}$,
$w=\tfrac{25}{2}$, and
$w=\tfrac{29}{2}$ respectively.
For each form we find eigenvalues
divisible by $11$. For the forms $\eta^{13}(24\tau)$,
$\eta^{13}(24\tau) E_4(24\tau)$,
and $\eta^{13}(24\tau) E_4^2(24\tau)$ we find that $\lambda_{19}$
is divisible by $11$. For $\eta^{13}(24\tau) E_6(24\tau)$ we find
that the eigenvalue $\lambda_5$ is divisible by $11$.
From \eqn{Nrmod11} and \eqn{Heckedef} it follows that
\beqs
N(r,11,\tfrac{1}{24}(5^3 \cdot 11\cdot 19^3 n + 1)) \equiv 0
\pmod{11}
\eeqs
for each $r$ provided that $(n,5)=(n,19)=1$. We replace
$n$ by $5\cdot 19\cdot 24\cdot n + c$ where $(c,5)=(c,19)=1$
and $5^3\cdot11\cdot19^3c\equiv -1\pmod{24}$. The smallest such $c$
is $c=11$. This gives 
$$
N(r,11,5^4 \cdot 11\cdot 19^4 n + 4322599) \equiv 0
\pmod{11},
$$
which is \eqn{nicecong11a}. We claim that
\beq
\eta^{13}(24\tau) F(24\tau) \mid U_{11} \equiv 0\pmod{11},
\mylabel{eq:E13Fmod11}
\eeq
for $F(\tau)=E_4(\tau)$, $E_6(\tau)$ or $E_4^2(\tau)$.
We observe that
\beqs
\eta^{13}(24\tau) F(24\tau) \equiv \frac{ \Delta(24\tau) F(24\tau)}
                                   {\eta(24\cdot11\tau)}
\pmod{11}
\eeqs
so that
\beqs
\eta^{13}(24\tau) F(24\tau) \mid U_{11}
\equiv \frac{ \Delta(24\tau) F(24\tau)\mid U_{11}}
            {\eta(24\tau)}
\pmod{11}.
\eeqs
Let $F=E_4$ so that $\Delta(\tau) E_4(\tau) \in S_{16}(1)$.
We easily find that
\beqs
\Delta(\tau) E_4(\tau) \mid U_{11} \equiv
\Delta(\tau) E_4(\tau) \mid T(11) \equiv 0
\pmod{11},
\eeqs
just by checking that the first coefficient is divisible by $11$
since $\dim S_{16}=1$ and $\Delta(\tau) E_4(\tau) \mid T(11)\in S_{16}(1)$.
The proves \eqn{E13Fmod11} for the case $F=E_4$. The other two
cases are analogous. Now \eqn{E13Fmod11} together with the fact that the 
eigenvalue $\lambda_{19}$ for $\eta^{13}(24\tau)$ is divisible by $11$
implies that
\beqs
N(r,11,\tfrac{1}{24}(11^2\cdot 19^3 n + 1)) \equiv 0
\pmod{11}
\eeqs
for each $r$ provided that $(n,19)=1$. We replace
$n$ by $19\cdot 24\cdot n + c$ where $(c,19)=1$
and $11^2\cdot19^3c\equiv -1\pmod{24}$. The smallest such $c$
is $c=5$. This gives
$$
N(r,11,11^2\cdot 19^4 n + 172904) \equiv 0
\pmod{11},
$$
which is \eqn{nicecong11b}.

\section{Explicit Congruences for $\spt(n)$} \label{sec:explicitsptcongs}

For $\ell > 3$ prime and $\beta \ge 1$ we define
\beqs
\SPT(\ell,\beta) = \sum_{n=0}^\infty \spt(\ell n + \beta) q^n.
\eeqs
We obtain a number of explicit congruences when $\ell \le 37$ and
$24\beta\equiv1\pmod{\ell}$.

\begin{theorem}
\mylabel{thm:SPTcongs}
We have
\begin{align}
\SPT(5,4) &\equiv 0 \pmod{5},\mylabel{eq:SPT5}\\
\SPT(7,5) &\equiv 0 \pmod{7},\mylabel{eq:SPT7}\\
\SPT(11,6) &\equiv 4 E^{13}(q) \pmod{11},\mylabel{eq:SPT11}\\
\SPT(13,6) &\equiv 0 \pmod{13},\mylabel{eq:SPT13}\\
\SPT(17,5) &\equiv 14 E^{7}(q) E_6(\tau) \pmod{17},\mylabel{eq:SPT17}\\
\SPT(19,4) &\equiv 10 E^{5}(q) E_4^2(\tau) \pmod{19},\mylabel{eq:SPT19}\\
\SPT(23,1) &\equiv E(q) E_4^3(\tau) + 7 q E^{25}(q) \pmod{23},\mylabel{eq:SPT23}\\
\SPT(29,23) &\equiv 17 E^{19}(q) E_6(\tau) \pmod{29},\mylabel{eq:SPT29}\\
\SPT(31,22) &\equiv 30 E^{17}(q) E_4^2(\tau) \pmod{31},\mylabel{eq:SPT31}\\
\SPT(37,17) &\equiv 12 E^{11}(q) E_4^2(\tau) E_6(\tau) \pmod{37},\mylabel{eq:SPT37}
\end{align}
\end{theorem}

Using Theorem \thm{SPTcongs},
Proposition \propo{Hecke} and Corollary \corol{Heckecor} we are able to
derive a number 
of explicit congruences for $\spt(n)$.                                   

\begin{theorem}
\mylabel{thm:sptexpcongs}
We have
\begin{align}
\spt( 11\cdot19^4\cdot n + 22006) &\equiv 0 \pmod{11}, \mylabel{eq:sptmod11}\\
\spt( 7^4\cdot17\cdot n + 243) &\equiv 0 \pmod{17}, \mylabel{eq:sptmod17}\\
\spt( 5^4\cdot19\cdot n + 99) &\equiv 0 \pmod{19}, \mylabel{eq:sptmod19}\\
\spt( 13^4\cdot29\cdot n + 18583) &\equiv 0 \pmod{29}, \mylabel{eq:sptmod29}\\
\spt( 29^4\cdot31\cdot n + 409532) &\equiv 0 \pmod{31}, \mylabel{eq:sptmod31}\\
\spt( 5^4\cdot37\cdot n + 1349) &\equiv 0 \pmod{37}. \mylabel{eq:sptmod37}
\end{align}
\end{theorem}

We prove Theorems \thm{SPTcongs} and \thm{sptexpcongs} case by case.
As usual we let $\ell > 3$ be prime 
and define $1\le \beta_{\ell}\le \ell-1$
such that $24\beta_{\ell}\equiv1\pmod{\ell}$.

We need some results from \cite{AtG} and \cite{BGM}.
Let $n$ be a positive integer. Define
\beqs
\mathcal{W}_{2n} =
\Span\{ \Phi_1^a \Phi_3^b \Phi_5^c \,:\, 1 \le a + 2b + 3c \le n \quad
\mbox{with $a$, $b$, $c$ nonnegative integers}\}
\eeqs
so that $\mathcal{W}_{2n}$ a vector space of quasimodular forms of bounded
weight over $\Q$. 
We now define by what we exactly mean by an  
``$\ell$-integral quasimodular form.''
We consider functions $E_{2}^{a}(z) \, F_b(z)$, where $F_b(z)\in M_b(1)$,
the coefficients in the $q$-expansion 
of $F_b(z)$ are $\ell$-integral and have bounded denominators, and
$a$ and $b$ are nonnegative integers. We call such a function
an $\ell$-integral quasi-modular form of weight $2a+b$.
Let $k$ be a nonnegative integer. In general, an  
\textit{$\ell$-integral quasi-modular form of weight $k$} is sum of such 
functions
where $2a+b=k$.
We let $\mathcal{X}_{2n}=\mathcal{X}_{2n,\ell}$ denote the subset 
of $\ell$-integral quasimodular forms in $\mathcal{W}_{2n}$.
By Theorems 7.4 and 7.6 in \cite{BGM} 
\beq
R_{2k} - P_{k}(\delta_q)R_2 \in P \mathcal{X}_{2k,\ell}
\subset P\mathcal{W}_{2k},
\mylabel{eq:R2Note1}
\eeq
where $P_{k}(x)\in \Z[x]$, for $2\le k \le \frac{\ell-3}{2}$ and $k=\frac{\ell+1}{2}$.
The result for $k=\frac{\ell-1}{2}$ is analogous except a factor $\ell$ may
occur in some denominators. 
We need a basis for $P\mathcal{W}_{2k}$ in terms of crank moments and
cusp forms. Following \cite[(5.14)]{AtG} we let
\beqs
\mathcal{C}_{2k} = \{ \delta_q^m(C_{2j})\,:\, 1\le j\le k,\, j+m\le k\}
\subset P\,\mathcal{W}_{2k}.
\eeqs
We let $\mathcal{B}_{2k}$ be a basis for the space of cusp forms $S_{2k}(1)$,
and let
\beqs
\mathcal{S}_{2k} = \cup_{\substack{1\le j\le k,\\ 0\le m\le k-j}}
                                    \delta_q^m(P\,\mathcal{B}_{2j})
\subset P\,\mathcal{W}_{2k}.
\eeqs
We
\begin{conj}
For $k\ge1$ the set
\beqs
\mathcal{T}_{2k}=\mathcal{C}_{2k} \cup \mathcal{S}_{2k}
\eeqs
forms a basis for $P\,\mathcal{W}_{2k}$ over $\Q$.
\end{conj}
We have confirmed this conjecture for $k\le 20$ which includes the cases
we need. We also note that 
\beqs
\dim (P\,\mathcal{W}_{2k}) = \abs{\mathcal{T}_{2k}},
\eeqs
using \cite[(3.32)]{AtG} and the fact that
\beqs
\dim M_{2k}(1) = 1 + \dim S_{2k}(1).
\eeqs

For $\ell>3$ prime and  $\epsilon\in\{-1,0,1\}$ we define the operator
$U_{\epsilon,\ell}^{*}$
which acts on $q$-series by
\beq
U_{\epsilon,\ell}^{*}\left(\sum a(n) q^n \right)
= \sum_{\leg{1-24n}{\ell}=\epsilon} a(n) q^n.
\mylabel{eq:Usdef}
\eeq
In \cite[Corollary 7.5]{BGM} 
we gave an elementary proof that for $\ell>3$ prime and  
$\epsilon=-1$ or $0$ there is a $G_{\ell}\in\mathcal{X}_{\ell+1,\ell}$
such that
\beq
U_{\epsilon,\ell}^{*}\left( {R}_2 \right)
\equiv
U_{\epsilon,\ell}^{*}\left( G_{\ell} P\right) \pmod{\ell}.
\mylabel{eq:bigR2cor}
\eeq
The proof depends on \eqn{R2Note1} and the fact that 
$R_{\ell+1}\equiv R_2\pmod{\ell}$, as well some elementary properties
of the polynomial $P_{\tfrac{\ell+1}{2}}(x)$.

\subsection*{Obtaining a congruence mod $\ell$ for $\SPT(\ell,\beta_\ell)$}
To prove Theorem \thm{SPTcongs} for the prime $\ell$ we follow
a series of computational steps:

\medskip

\noindent
\textit{Step 1}. 
First we compute 
\beqs
L_{2k} := R_{2k} - P_{k}(\delta_q)R_2 
\eeqs
in terms of the set $\mathcal{T}_{2k}$ for $k=\tfrac{\ell-1}{2}$ and 
$k=\tfrac{\ell+1}{2}$. This step involves the heaviest computation.
In each case this calculation was done in
MAPLE by computing the coefficients of $q^j$ ($0\le j\le n$ where 
$n=\abs{\mathcal{T}_{2k}}+20$) of each function in 
$\mathcal{T}_{2k}$ as well as the function $L_{2k}$. These
coefficients form a $(n+20)\times n$ matrix $A$, where each
column corresponds to a function in $\mathcal{T}_{2k}$ with the last column
corresponding to the function $L_{2k}$.      We used MAPLE to show that
$\dim \mbox{Nul}(A)=1$ and find a basis vector for this Nullspace.
In each case the last component of 
this basis vector is nonzero thus giving  the function $L_{2k}$ as a 
linear combination
of the functions in $\mathcal{T}_{2k}$ since we know that
$L_{2k}\in P\mathcal{W}_{2k}$ and $\abs{\mathcal{T}_{2k}}=\dim P\mathcal{W}_{2k}$. 

\medskip

\noindent
\textit{Step 2}. 
For the identity corresponding to $k=\tfrac{\ell-1}{2}$ we find coefficients
are $\ell$-integral except for a factor $\ell$ in some denominators.
We multiply both sides of the identity by $\ell$, apply the $U_{0,\ell}^{*}$
operator (see \eqn{Usdef}) and reduce mod $\ell$. This gives an identity 
mod $\ell$ for
the second crank moment mod $\ell$ in terms of certain cusp forms times $P$.
Alternatively, this can be computed by using \eqn{DyM2}.

\medskip

\noindent
\textit{Step 3}. 
For the identity corresponding to $k=\tfrac{\ell+1}{2}$ we apply the $U_{0,\ell}^{*}$
operator and reduce mod $\ell$. 
This gives an identity mod $\ell$ for
the second rank moment mod $\ell$ in terms of certain cusp forms times $P$.
See \eqn{bigR2cor}.

\medskip

\noindent
\textit{Step 4}. 
Using \eqn{sptid} and the identities from Steps 2 and 3 we obtain a
congruence mod ${\ell}$ for $\SPT(\ell,\beta_\ell)$ in terms of 
$U_{0,\ell}^{*}(P\,F)$ 
where $F$ is a sum of cusp forms of different weights. We use the theory of 
modular forms mod $\ell$ and Theorem \thm{pFcong} to simplify and obtain the 
final result.  In each case we find that 
\beqs
\SPT(\ell,\beta_\ell) \equiv E^{r_\ell}(q) F_\ell(\tau) \pmod{\ell}
\eeqs
for some $F_\ell(\tau)\in M_{\tfrac{1}{2}(\ell-r_\ell)+1}(1) \cap \Z[[q]]$.

\medskip

\subsection*{Obtaining an explicit  congruence mod $\ell$ for $\spt(n)$} 
We let 
\beqs
K_{\ell}(\tau) = \eta^{r_\ell}(24\tau)  F_\ell(24\tau) \in 
S_{\tfrac{\ell+2}{2}}(576,\chi_{12}),
\eeqs
$F_\ell(\tau)$ was found in Step 4 (above). In most cases the
function $K_{\ell}(\tau)$ is a Hecke eigenform.

\medskip

\noindent
\textit{Step 5}. Suppose  $K_{\ell}(\tau)$ is a Hecke eigenform; i.e.
for each prime $Q>3$
\beqs
K_\ell(\tau)\mid T(Q^2) = \lambda_Q K_\ell(\tau),
\eeqs
for some integer $\lambda_Q$.
Find the first prime $Q$ such the eigenvalue $\lambda_Q\equiv0\pmod{\ell}$.
Now
\beqs
\sum_{n=0}^\infty \spt(\ell n + \beta_\ell) q^{24n+r_\ell}
= \sum_{n=0}^\infty \spt\left(\frac{\ell n + 1}{24}\right) q^n
\equiv \eta^{r_\ell}(24\tau)  F_\ell(24\tau) \equiv K_\ell(\tau)
\pmod{\ell}.
\eeqs
It follows from \eqn{Heckedef} that
\beq
\spt\left(\frac{\ell Q^3 n + 1}{24}\right) \equiv 0 \pmod{\ell}
\mylabel{eq:Step5B}
\eeq
provided that $(n,Q)=1$. 

\medskip

\noindent
\textit{Step 6}. For the prime $Q$ found in Step 5 find the smallest
integer $c$, such that $(c,Q)=1$ and $c\equiv -\ell Q\pmod{24}$.
Then in \eqn{Step5B} we replace $n$ by $24Qn +c$ to obtain
\beq
\spt(\ell Q^4 n + \tfrac{1}{24}(\ell Q^3 c + 1)) \equiv 0 \pmod{\ell}.
\mylabel{eq:expcong}
\eeq
We note that the constant $\tfrac{1}{24}(\ell Q^3 c + 1)$ is an
integer since $\ell^2\equiv Q^2\equiv 1\pmod{24}$ and
$c\equiv -\ell Q\pmod{24}$.

We now prove  Theorems \thm{SPTcongs} and \thm{sptexpcongs} case by case by 
following Steps 1--6. All calculations were done using MAPLE.
In Section \sect{spt11} we gave a proof of the $\ell = 11$ case.
Here we give the proof of the $\ell = 17$ case in detail and sketch the
other cases.

\subsubsection*{$\ell = 17$}
We compute
$$
L_{16}:=R_{16} - P_{8}(\delta_q)R_2
$$
in terms of the $40$ functions in $\mathcal{T}_{16}$.Then we 
multiply by $17$, apply $U_{0,17}^{*}$ and reduce mod $17$ to find
\beq
U_{0,17}^{*}\left(C_2\right) \equiv
U_{0,17}^{*}\left( 8 E_4 \Delta P\right)
\pmod{17}.
\mylabel{eq:R16mod17}
\eeq
Next we compute
$$
L_{18}:=R_{18} - P_{9}(\delta_q)R_2
$$
in terms of the $52$ functions in $\mathcal{T}_{18}$. Then 
apply $U_{0,17}^{*}$, reduce mod $17$ and use \eqn{R16mod17} to find
\beqs
\Uop{17}{R_{2}} \equiv \Uop{17}{
 8  E_4 \Delta P + 11 E_6 \Delta P}
\pmod{17}.
\eeqs
Thus we have
\beqs
\SPT(17,5) = \tfrac{1}{2} \Uop{17}{C_2 - R_2}
\equiv \Uop{17}{3 E_6 \Delta P}.
\eeqs
By Theorem \thm{pFcong}
\beqs
\Uop{17}{3 E_6 \Delta P} \equiv E^7(q) G(\tau) \pmod{17},
\eeqs
for some $G(\tau)\in M_{22}(1)\cap\Z[[q]]$. A finite 
computation shows that
\beq
G(\tau) \equiv 14 E_{16}(\tau) E_6(\tau) \pmod{17}.
\mylabel{eq:Gmod17}
\eeq 
In fact we need only verify the coefficients of $q^n$ 
on both sides of \eqn{Gmod17} agree mod $17$ for
$n\le \lfloor \tfrac{22}{12}\rfloor = 2$. This is a standard argument.
See for example \cite{St}.
But by \cite[Theorem 2(i)]{SWD} 
\beqs
E_{16}(\tau) \equiv 1 \pmod{17},
\eeqs
and we have
\beqs
\SPT(17,5) \equiv 14 E^7(q) E_{16}(\tau) E_6(\tau) 
\equiv 14 E^7(q) E_6(\tau) \pmod{17},
\eeqs
which is \eqn{SPT17}.
Thus
\beqs
\sum_{n=0}^\infty \spt(17 n + 5) q^{24n+7}
\equiv 14 G_{17}(\tau)
\pmod{17},
\eeqs
where
\beqs
G_{17}(\tau) = \eta^7(24\tau) E_6(24\tau)\in S_{\tfrac{19}{2}}(576,\chi_{12}). 
\eeqs
Since $\dim M_6(1)=1$, $G_{17}(\tau)$ is a Hecke eigenform by Corollary
\corol{Heckecor}. The first eigenvalue divisible by $17$ is
\beqs
\lambda_7 = -24959264 = - 2^5 \cdot 11 \cdot 17 \cdot 43 \cdot 97
\eeqs
We want the smallest integer $c$ satisfying $(c,7)=1$ and
$c\equiv (-17)7 \equiv 1\pmod{24}$; i.e.\ $c=1$.
As in \eqn{expcong} we obtain the congruence
\beqs
\spt( 7^4\cdot 17 \cdot n + 243) \equiv 0 \pmod{17},
\eeqs
which is \eqn{sptmod17}.

\subsubsection*{$\ell = 19$}
We find
\begin{align*}
U_{0,19}^{*}\left(C_2\right) &\equiv
U_{0,19}^{*}\left( 5 E_6 \Delta P\right)
\pmod{19},  \\
\Uop{19}{R_{2}} &\equiv \Uop{19}{
 5  E_6 \Delta P +  E_4^2 \Delta P}
\pmod{19}, \\
\SPT(19,4) &= \tfrac{1}{2} \Uop{19}{C_2 - R_2}
\equiv \Uop{19}{9 E_4^2 \Delta P}.
\end{align*}
Using Theorem \thm{pFcong} we find that
\beqs
\Uop{19}{9 E_4^2 \Delta P} \equiv 
10 E^5(q) E_{18}(\tau) E_4^2(\tau) \equiv 10 E^5(q) E_4^2(\tau) \pmod{19}.
\eeqs
The result \eqn{SPT19} follows.
Thus
\beqs
\sum_{n=0}^\infty \spt(19 n + 4) q^{24n+5}
\equiv 10 G_{19}(\tau)
\pmod{19},
\eeqs
where
\beqs
G_{19}(\tau) = \eta^5(24\tau) E_4^2(24\tau)\in S_{\tfrac{21}{2}}(576,\chi_{12}).
\eeqs
Since $\dim M_8(1)=1$, $G_{19}(\tau)$ is a Hecke eigenform by Corollary
\corol{Heckecor}. The first eigenvalue divisible by $19$ is $\lambda_5$
and \eqn{sptmod19} follows.

\subsubsection*{$\ell = 23$}
We find
\begin{align*}
U_{0,23}^{*}\left(C_2\right) &\equiv
U_{0,23}^{*}\left(2 E_{10} \Delta P\right)
\pmod{23},  \\
\Uop{23}{R_{2}} &\equiv \Uop{23}{
 2  E_{10} \Delta P +  9 \Delta^2 P + 21 E_4^3 \Delta P}
\pmod{23}, \\
\SPT(23,1) &= \tfrac{1}{2} \Uop{23}{C_2 - R_2}
\equiv \Uop{23}{7\Delta^2 P + E_4^3 \Delta P}.
\end{align*}
By Theorem \thm{pFcong}
\beqs
\Uop{23}{7\Delta^2 P + E_4^3 \Delta P} 
\equiv E(q) E_{22}(\tau) (7\Delta(\tau) + E_4^3(\tau))
        \equiv E(q) (7\Delta(\tau) + E_4^3(\tau))
 \pmod{23}.
\eeqs
The result \eqn{SPT23}
follows. However in this case the function 
$E(q) (7\Delta(\tau) + E_4^3(\tau))$ is not a Hecke eigenform.

\subsubsection*{$\ell = 29$}
We find
\begin{align*}
U_{0,29}^{*}\left(C_2\right) &\equiv
U_{0,29}^{*}\left(\Delta P(2 E_4^2 + 20 E_4^4 + 11 \Delta E_4)\right)
\pmod{29}, \\
\Uop{29}{R_{2}} &\equiv \Uop{29}{
\Delta P(25 E_4^2 + 20 E_4^4 + 11 \Delta E_4 + 5 \Delta E_6)}
\pmod{29}, \\
\SPT(29,23) &= \tfrac{1}{2} \Uop{29}{C_2 - R_2}
\equiv \Uop{29}{\Delta P( 3 E_4^2 + 12 \Delta E_6)}
\pmod{29}.
\end{align*}
Using Theorem \thm{pFcong} we find
$\Uop{29}{E_4^2 \Delta P} \equiv 0 \pmod{29}$ 
and
\beqs
\Uop{29}{12 E_6 \Delta^2 P} 
\equiv 17 E^{19}(q) E_{28}(\tau) E_6(\tau) 
\equiv 17 E^{19}(q) E_6(\tau) \pmod{29}.
\eeqs
The result \eqn{SPT29}
follows.
Thus
\beqs
\sum_{n=0}^\infty \spt(29 n + 23) q^{24n+23}
\equiv 17 G_{29}(\tau)
\pmod{29},
\eeqs
where
\beqs
G_{29}(\tau) = \eta^{19}(24\tau) E_6(24\tau)\in 
S_{\tfrac{31}{2}}(576,\chi_{12}).
\eeqs
Since $\dim M_6(1)=1$, $G_{29}(\tau)$ is a Hecke eigenform by Corollary
\corol{Heckecor}. The first eigenvalue divisible by $29$ is $\lambda_{13}$
and \eqn{sptmod29} follows.

\subsubsection*{$\ell = 31$}
We find
\begin{align*}
U_{0,31}^{*}\left(C_2\right) &\equiv
U_{0,31}^{*}\left(
2 \Delta P \left( 3 E_{{6}}\Delta + 6 {E_{{6}}}^{3}
                 +2 E_{{6}}+5 E_{{10}} \right) \right)
\pmod{31}, \\
U_{0,31}^{*}\left(R_2\right) &\equiv
\Delta P \left( 6 E_{{6}}\Delta+29 \Delta {E_{{4}}}^{2}+
                12 {E_{{6}}}^{3}+2 E_{{6}}+30 E_{{10}} \right) 
\pmod{31}, \\
\SPT(31,22) &= \tfrac{1}{2} \Uop{31}{C_2 - R_2}
\equiv \Uop{31}{\Delta P \left( \Delta {E_{{4}}}^{2}
                  +E_{{6}}+21 E_{{10}} \right)}
\pmod{31}.
\end{align*}
By using 
Theorem \thm{pFcong} we show 
that
$\Uop{31}{E_6    \Delta P} \equiv 
\Uop{31}{E_{10} \Delta P} \equiv 0 \pmod{31}$,
and 
\beqs
\Uop{31}{E_4^2 \Delta^2 P} 
\equiv 30 E^{17}(q) E_{30}(\tau) E_4^2(\tau) 
\equiv 30 E^{17}(q) E_4^2(\tau) \pmod{31}.
\eeqs
The result \eqn{SPT31}
follows.
Thus
\beqs
\sum_{n=0}^\infty \spt(31 n + 22) q^{24n+17}
\equiv 30 G_{31}(\tau)
\pmod{31},
\eeqs
where
\beqs
G_{31}(\tau) = \eta^{17}(24\tau) E_4^2(24\tau)\in
S_{\tfrac{33}{2}}(576,\chi_{12}).
\eeqs
Since $\dim M_8(1)=1$, $G_{31}(\tau)$ is a Hecke eigenform by Corollary
\corol{Heckecor}. The first eigenvalue divisible by $31$ is
$\lambda_{29}$
and \eqn{sptmod31} follows.

\subsubsection*{$\ell = 37$}
We find    
\begin{align*}
U_{0,37}^{*}\left(C_2\right) &\equiv
U_{0,37}^{*}\left(
\Delta P \left( 22 {\Delta}^{2}+4 \Delta+31 \Delta {E_{{6}}}^{2}
        +17 \Delta E_{{4}}+4 {E_{{4}}}^{3}+3 E_{{10}} \right.\right.
 \\
       &\quad  \left.\left. +27 {E_{{4}}}^{4}
        +12 {E_{{4}}}^{6} \right) \right)
\pmod{37},\\
\Uop{37}{R_{2}} &\equiv 
U_{0,37}^{*}\left(
\Delta P \left( 22 {\Delta}^{2}+36 \Delta+31 \Delta {E_{{6}}}^{2}
+35 \Delta E_{{4}}+20 \Delta E_{{6}}{E_{{4}}}^{2}+36 {E_{{4}}}^{4}
\right.\right.
 \\
&\left.\left.\quad+13 E_{{14}}{E_{{4}}}^{3}+23 E_{{10}}+36 {E_{{4}}}^{3}
+12 {E_{{4}}}^{6} \right) \right)
\pmod{37},
\\
\SPT(37,17) &= \tfrac{1}{2} \Uop{37}{C_2 - R_2}  \\
&\equiv 
U_{0,37}^{*}\left(
\Delta P \left( 21 \Delta+28 \Delta E_{{4}}+27 \Delta E_{{6}}{E_{{4}}}^{2}
+27 E_{{10}}+14 {E_{{4}}}^{4}+21 {E_{{4}}}^{3}+
12 E_{{14}}{E_{{4}}}^{3} \right) \right)
\pmod{37}.
\end{align*}
Using Theorem \thm{pFcong} we find that
\beqs             
\Uop{37}{E_{10}\Delta P} \equiv
\Uop{37}{\Delta P( \Delta + E_4^3)} \equiv
\Uop{37}{\Delta P( 28 \Delta E_4 + 14 E_4^4)} 
\equiv 0 \pmod{37},
\eeqs
and                                         
\beqs
\Uop{37}{\Delta P( 27\Delta E_6 E_4^2 + 12 E_{14} E_4^3)} 
\equiv  12 E^{11}(q) E_{36}(\tau) E_4^2(\tau) E_6(\tau)
\equiv 12 E^{11}(q) E_4^2(\tau) E_6(\tau)
\pmod{37}.
\eeqs
The result \eqn{SPT37}
follows.
Thus
\beqs
\sum_{n=0}^\infty \spt(37 n + 17) q^{24n+11}
\equiv  12 G_{37}(\tau)
\pmod{37},
\eeqs
where
\beqs
G_{37}(\tau) = \eta^{11}(24\tau) E_4^2(24\tau) E_6(24\tau) \in
S_{\tfrac{39}{2}}(576,\chi_{12}).
\eeqs
Since $\dim M_{14}(1)=1$, $G_{37}(\tau)$ is a Hecke eigenform by Corollary
\corol{Heckecor}. The first eigenvalue divisible by $37$ is $\lambda_5$
and \eqn{sptmod37} follows.
\section{Concluding remarks} \label{sec:conclude}
In a companion paper \cite{BGM} we further explore the general problem of 
congruences mod $\ell$ for rank moments, Andrews \cite{Adurf} symmetrized 
rank moments and full rank
functions for $k$-marked Durfee symbols for general $\ell$.
Let $\ell>3$ be prime, and $\beta_\ell$, $r_\ell$ be as in
Theorem \thm{pFcong}. By a detailed $\ell$-adic analysis of
the rank-crank moment equation \eqn{rcrel} we have been
able to show that
\beq
\sum_{n=0}^\infty N_{2k}(\ell n + \beta_\ell) q^{24n + r_\ell}
\equiv \eta^{r_\ell}(24\tau)\, G_{\ell, 2k}(24\tau)
\pmod{\ell},
\mylabel{eq:N2kgencong}
\eeq
where $G_{\ell,2k}(\tau)$ is a sum of level one integral modular forms
of bounded weight when $2 \le 2k \le \ell -3$. When $2k=\ell -1$ the result
is similar but involves an additional function analogous to $\Psi$, given
in \eqn{Psidef}. By using a result of Treneer \cite[Prop.~5.2]{Tren}
it can be shown that
there exists a positive proportion of the primes $Q$ for which the appropriate
Hecke operator $T(Q^2)$ annihlates mod $\ell$ all the half integer weight
forms that occur on the right side of \eqn{N2kgencong}. This implies that
there exists a positive proportion of the primes $Q$ 
such that
\beq
N_{2k}(\ell Q^4 n +  \mu) \equiv 0 \pmod{\ell},
\mylabel{eq:NKcongs}
\eeq
for all integers $k$, $n\ge0$. Here $\mu$ is an integer constant
that depends on $\ell$ and $Q$. 
As a result this leads to infinitely many Ramanujan-type congruences
for Dyson's rank function $N(r,\ell,n)$ (the number of partitions
of n with rank congruent to $r$ mod $\ell$).
In particular,  this implies that for each prime
$\ell>3$ there is a positive proportion of the primes $Q$
such that
\beq
N(r,\ell,\ell Q^4 n +  \mu) \equiv 0 \pmod{\ell},
\mylabel{eq:rankcongsIntro}
\eeq
for all integers $n\ge0$ and $0\le r\le \ell-1$.
Again $\mu$ is an integer constant
at depends on $\ell$ and $Q$.  This corresponds to Bringmann's \cite{Bringmann}
result \eqn{NrAnBcongs} (with $m=u=1$) which was proved using the theory
of Maass forms.
As noted before, since $\spt(n)$ involves $p(n)$ and the second
rank moment $N_2(n)$ this leads to another proof  that there are
infinitely many Ramanujan-type congruences mod $\ell$ ($\ell>3$
prime) for $\spt(n)$.   
In particular,  this implies that for each prime
$\ell>3$ there is a positive proportion of the primes $Q$
such that
\beq
\spt(\ell Q^4 n +  \mu) \equiv 0 \pmod{\ell},
\mylabel{eq:sptcongsexplicit}
\eeq
for all integers $k$, $n\ge0$ and $\mu$ is an integer constant
at depends on $\ell$ and $Q$. We note that Bringmann \cite{Bringmann2}
proved this result by a different method.

In a recent paper, Bringmann, Ono and Rhoades \cite{BOR} have shown the
generating function for many rank differences is a modular form.
In particular, they have shown the function
\beq
\sum_{n=0}^\infty \left(N(r_1,\ell,\ell n + d) -
N(r_2,\ell,\ell n + d) \right) q^{24(\ell n + d) -1}
\mylabel{eq:rankdiffs}
\eeq
is a weakly holomorphic weight $1/2$ modular form when
$1-24d$ is a quadratic nonresidue mod $\ell$ and when it is
a quadratic residue under certain additional conditions on $r_1$
and $r_2$.
As noted in \cite{BOR} this was observed earlier by Atkin and
Swinnerton-Dyer \cite{ASD} in the cases $\ell=5$ and $\ell=7$.
It should also be noted that it also been observed when
$\ell=11$ by Atkin and Hussain \cite{AH}, and $\ell=13$
by O'Brien \cite{OB}. What is more interesting is the result
also holds in the case $24d\equiv1\pmod{\ell}$ for $\ell=11$ 
and $\ell=13$ (again \cite{AH} and \cite{OB}). We conjecture that the analogue
for the case $24d\equiv1\pmod{\ell}$
holds for all primes $\ell > 3$. The cases $\ell=5$, $7$ are trivial
since the rank differences are zero. We can prove that in the case
$24d\equiv1\pmod{\ell}$ the generating function for the rank differences
is congruent mod $\ell$ to a sum of half-integer weight cusps. 
An analog of this also holds in the case that $1-24d$ is a quadratic 
nonresidue mod $\ell$.

\noindent
\textbf{Acknowledgement}

\noindent
I would like to thank the following for helpful comments: 
Kathrin Bringmann, Byungchan Kim, Karl Mahlburg, and Ken Ono.


\bibliographystyle{amsplain}

\begin{thebibliography}{10}
\bibitem{AB}
S.~Ahlgren and M.~Boylan, 
\textit{Arithmetic properties of the partition function},
Invent. Math.
\textbf{153} 
(2003),
487--502.
\bibitem{Adurf}
G.~E.~Andrews,
\textit{Partitions, Durfee symbols, and the Atkin-Garvan moments of ranks},
Invent. Math. 
\textbf{169}
(2007), 
37--73.                                 
\bibitem{Aspt}
G.~E.~Andrews,
\textit{The number of smallest parts in the partitions of $n$},
J.\ Reine Angew.\ Math., to appear.
\bibitem{AG}
G.~E.~Andrews and F.~G.~Garvan,
\textit{Dyson's crank of a partition},
{Bull. Amer. Math. Soc. (N.S.)}
\textbf{18} (1988), 167--171.
%
\bibitem{Ap}
T.~M.~Apostol,
\textit{Modular Functions and Dirichlet Series in Number Theory},
Springer-Verlag, New York, 1976.
\bibitem{At}
A.~O.~L.  Atkin, 
\textit{Multiplicative congruence properties and density problems for
$p(n)$},
Proc. London Math. Soc. (3)
\textbf{18} 
(1968),
563--576.
\bibitem{At1996}
A.~O.~L.  Atkin, 
\textit{$p(n)$ Revisited}, NMBRTHY Archive, July 1996,
\texttt{http://listserv.nodak.edu/archives/nmbrthry.html}
\bibitem{AtG}
A.~O.~L.~Atkin and  F.~G.~Garvan,         
\textit{Relations between the ranks and cranks of partitions},
{Ramanujan J.}
\textbf{7} (2003), {343--366}.
\bibitem{AH}
A.~O.~L.~Atkin and S.~M.~Hussain,
\textit{Some properties of partitions II},
Trans. Amer. Math. Soc.  \textbf{89} (1958), 184--200.
\bibitem{AtkinLehner}
A.~O.~L.~Atkin, and J.~Lehner, 
\textit{Hecke operators on {$\Gamma \sb{0}(m)$}},
Math. Ann.
\textbf{185} 
(1970),
134--160.
%
\bibitem{ASD}
A.~O.~L.~Atkin and P.~Swinnerton-Dyer,
\textit{Some properties of partitions},
Proc. London Math. Soc.
\textbf{4} (1954), 84--106.
\bibitem{Bringmann}
K.~Bringmann, \textit{Congruences for Dyson's ranks}, preprint.
\textit{On certain congruences for Dyson's ranks},
Int. J. Number Theory, to appear.
\bibitem{Bringmann2}
K.~Bringmann, 
\textit{On the explicit construction of higher deformations of
       partition statistics},
       Duke Math. J., to appear.
preprint.                      
\bibitem{BGM}
K.~Bringmann, F.~G.~Garvan and K.~Mahlburg 
\textit{Partition statistics and quasiweak Maass forms}, preprint,
arXiv:0803.1891
\bibitem{BringmannOno}
{K.~Bringmann and K.~Ono},
\textit{Dyson's ranks and Maass forms},
Ann. of Math., to appear.
\bibitem{BOR}
{K.~Bringmann, K.~Ono and R.~C.~Rhoades},
\textit{Eulerian series as modular forms},
J. Amer. Math. Soc.,
to appear.
\bibitem{Bruinier}
J.~H.~Bruinier,
\textit{Nonvanishing modulo {$l$} of {F}ourier coefficients
of half-integral weight modular forms},
Duke Math. J.
\textbf{98}
(1999),
595--611.
\bibitem{Chua}
K.~S.~Chua, 
\textit{Explicit congruences for the partition function modulo every prime},
Arch. Math. (Basel)
\textbf{81} 
(2003),
11--21.
\bibitem{Dyson44}
F.~J.~Dyson,
\textit{Some guesses in the theory of partitions},
Eureka (Cambridge) \textbf{8} (1944), 10--15.
\bibitem{Dy}
F.~J.~Dyson, 
\textit{Mappings and symmetries of partitions},
J. Combin. Theory Ser. A
\textbf{51}
(1989),
169--180.
\bibitem{FolsomOno}
A.~Folsom and K.~Ono,
\textit{The spt-function of Andrews}, preprint.
\bibitem{G88a}
F.~G.~Garvan, 
\textit{New combinatorial interpretations of Ramanujan's partition 
congruences mod $5,7$ and $11$}, 
Trans. Amer. Math. Soc.  \textbf{305}  (1988), 47--77. 
\bibitem{Mahlburg}
K.~Mahlburg, 
\textit{Partition congruences and the {A}ndrews-{G}arvan-{D}yson crank},
Proc. Natl. Acad. Sci. USA
\textbf{102} 
(2005),
15373--15376.
\bibitem{OB}
 J.~N.~O'Brien,
\textit{Some properties of partitions with special reference to primes
other than $5$, $7$ and $11$},
Ph.D. thesis, Univ. of Durham, England, 1966.
\bibitem{Ono}
K.~Ono, 
\textit{Distribution of the partition function modulo {$m$}}
Ann. of Math. (2)
\textbf{151} 
(2000),
293--307.
\bibitem{Guo-Ono}
L. Guo and K. Ono,
\textit{The partition function and the arithmetic of certain modular
$L$-functions}
Internat. Math. Res. Notices, 1999, No. 21,
1179--1197.
\bibitem{Ram}
S.~Ramanujan,
\textit{On certain arithmetic functions}
Trans Cambridge Philos. Soc.
\textbf{XXII} (1916), 159--184.
\bibitem{SWD}
{H.~P.~F.~Swinnerton-Dyer}, 
\textit{On {$l$}-adic representations and congruences for coefficients
              of modular forms}, in
``Modular functions of one variable, III'' (Proc. Internat. Summer
              School, Univ. Antwerp, 1972),
Springer, Berlin, 1973, pp. 1--55.
\bibitem{Sh}
G.~Shimura, 
\textit{On modular forms of half integral weight}
Ann. of Math. (2)
\textbf{97}
(1973),
440--481.
\bibitem{St}
J.~Sturm,
\textit{On the congruence of modular forms},
in ``Number Theory,'' 
{Lecture Notes in Math.},
{vol. 1240},
{Springer},
{Berlin},
{1987}, pp. 275--280.
\bibitem{Tren}
S. Treneer 
\textit{Quadratic twists and the coefficients of weakly holomorphic
modular forms}, preprint.
\bibitem{We}
R.~L.~Weaver,
\textit{New congruences for the partition function}
Ramanujan J.
\textbf{5} 
(2001),
53--63.




\end{thebibliography}

\end{document}